\documentclass[10pt,a4paper]{article}
\usepackage{amsfonts}
\usepackage{latexsym}
\usepackage{cite}
\usepackage{amsmath,amsfonts,latexsym,amssymb}
\usepackage[mathscr]{eucal}
\usepackage{cases,color}
\usepackage{amsthm}

\usepackage[bf,small]{caption2}
\usepackage{float}
\usepackage{graphicx}
\usepackage{amsmath}
\usepackage{amssymb}
\usepackage[all]{xy}

\newtheorem{theorem}{theorem}[section]

\newtheorem{lem}[theorem]{Lemma}

\newtheorem{que}[theorem]{Question}
\newtheorem{rmk}[theorem]{Remark}
\newtheorem{thm}[theorem]{Theorem}
\newtheorem{ass}[theorem]{Assertion}

\begin{document}

\title{\vspace{-2cm}\textbf{The ${\rm SL}(2,\mathbb{C})$-character variety of $8_{18}$}}
\author{\Large Haimiao Chen}

\date{}
\maketitle

\begin{abstract}
  The knot $8_{18}$ is the first non-arborescent hyperbolic knot. In 2020, Paoluzzi and Porti found its ${\rm SL}(2,\mathbb{C})$-character variety with the aid of a computer, but many details were omitted.
  In this paper, we determine the character variety by a software-free procedure, which is easy to follow and enlightening.
  Along the way, we develop an efficient method for working with simultaneous conjugacy classes of four elements of ${\rm SL}(2,\mathbb{C})$.

\end{abstract}

\section{Introduction}

Fix $G={\rm SL}(2,\mathbb{C})$. Suppose $\Gamma$ is a finitely presented group. Call an element $\rho\in\hom(\Gamma,G)$ a $G$-{\it representation} of $\Gamma$, and define its {\it character} as the function $\chi_\rho:\Gamma\to\mathbb{C}$ sending $x$ to ${\rm tr}(\rho(x))$.
Call $\rho$ reducible if the elements of its image have a common eigenvector, and irreducible otherwise.
As a well-known fact, two irreducible representations $\rho,\rho'$ are conjugate (i.e. there exists $\mathbf{a}\in G$ such that $\rho'(x)=\mathbf{a}\rho(x)\mathbf{a}^{-1}$ for all $x$) if and only if $\chi_\rho=\chi_{\rho'}$.
Call $\mathcal{X}(\Gamma):=\{\chi_\rho\colon\rho\in\hom(\Gamma,G)\}$ the character variety of $\Gamma$, and call the subset
$\mathcal{X}^{\rm irr}(\Gamma):=\{\chi_\rho\colon\rho\ \text{is\ irreducible}\}$ the {\it irreducible character variety}.

Character varieties of fundamental groups of 3-manifolds are prevalent in the literature on low-dimensional topology \cite{Sh96,Si12}.
Nowadays a consensus is that for a $3$-manifold $M$, much geometric and topological information on $M$ is carried by $\mathcal{X}(\pi_1(M))$.
Studying the character variety is a manner of understanding $M$.
For a knot $K\subset S^3$, let $E_K$ denote its exterior, and let $\pi(K)=\pi_1(E_K)$.
Let $\mathcal{X}(K)=\mathcal{X}(\pi(K))$, and $\mathcal{X}^{\rm irr}(K)=\mathcal{X}^{\rm irr}(\pi(K))$.
We focus on $\mathcal{X}^{\rm irr}(K)$, discarding reducible characters, as they are relatively easy to investigate.

\begin{figure}[h]
  \centering
  \includegraphics[width=2.8cm]{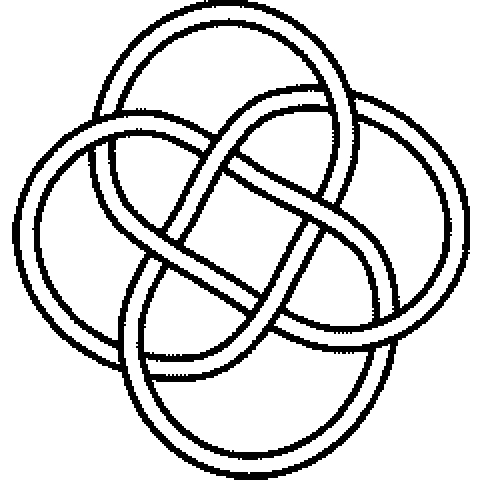}\\
  \caption{The knot $8_{18}$ (taken from Knot Atlas).}\label{fig:8_18}
\end{figure}

In this paper, we are interested in $8_{18}$, which is denoted by $T$ throughout.
It is a member of the family of Turk's head knots, denoted by $Th(3,4)$ in the notation of \cite{PS24}, and sometimes simply
called {\it the Turk's head knot} \cite{HRS12}.
It is the first non-arborescent knot in Rolfsen's table, and it is hyperbolic.

The hyperbolicity means that there exists a discrete faithful representation $\pi(T)\to{\rm PSL}(2,\mathbb{C})$, called the {\it holonomy representation}, which can be lifted to $\rho_\star:\pi(T)\to G$.
The irreducible component $\mathcal{X}_\star$ of $\mathcal{X}(T)$ containing $\chi_{\rho_\star}$ is called the {\it excellent component}.
It was found by Hilden, Lozano and Montesinos-Amilibia \cite{HLM00}, by making good use of the $4$-periodicity.

The whole character variety $\mathcal{X}(T)$ was determined in 2020 by Porti and Paoluzzi \cite{PP20}, with the aid of symbolic software, but many details were omitted.

We feel it reasonable to write a new paper, in several aspects.
An effective method of computing $\mathcal{X}^{\rm irr}(K)$ for arborescent knots $K$ had been developed in \cite{Ch23}, so we are curious about
non-arborescent knots. Since $T$ is a very interesting knot, it is worthwhile to thoroughly clarify the steps of finding $\mathcal{X}^{\rm irr}(T)$, through a procedure as elementary as possible.
While the approach of \cite{PP20} used $3$ generators for $\pi(T)$, we adopt a $4$-generator presentation which is compatible with the rotational symmetry, so as to describe $\mathcal{X}^{\rm irr}(T)$ as a symmetric affine algebraic set.
Beneficially, we are able to parameterize each component in a moderately succinct way, and present $\mathcal{X}^{\rm irr}(T)$ more explicitly.

More importantly, this work illustrates how to efficiently deal with simultaneous conjugacy classes of four $2\times 2$ matrices.
Let $F_n$ denote the free group of rank $n$. It is known that the coordinate ring $\mathbb{C}[\mathcal{X}(F_n)]$ is generated by
$n+{n\choose 2}+{n\choose 3}$ trace functions. The number of definitive relations for $\mathbb{C}[\mathcal{X}(F_3)]$ is $1$, but that for $\mathbb{C}[\mathcal{X}(F_4)]$ dramatically increases to 20. However, we point out that in most cases, actually it is enough to treat $5$ relations.
We develop new computational techniques, and expect them to be applicable in other occasions in future.

The content is organized as follows. In Section 2 we collect some preliminaries about $2\times 2$ matrices, and outline the strategy of computing $\mathcal{X}^{\rm irr}(T)$. In Section 3 we explicitly write down the trace equations expressing conditions on four given matrices to define an irreducible representation. In Section 4 we find all the solutions to these equations. Section 5 presents the main result (Theorem \ref{thm:main}) and some discussions.

\section{Preparation}

Let $\mathcal{M}$ denote the set of $2\times 2$ matrices with entries in $\mathbb{C}$; it is a $4$-dimensional vector space over $\mathbb{C}$.
Let $\mathcal{M}_0=\{\mathbf{x}\in\mathcal{M}\colon{\rm tr}(\mathbf{x})=0\}$.

For $\mathbf{x}\in\mathcal{M}$, denote its transpose and adjoint as $\mathbf{x}^{\rm tr}$ and $\mathbf{x}^\ast$, respectively.
Clearly, $\mathbf{x}+\mathbf{x}^\ast={\rm tr}(\mathbf{x})\cdot\mathbf{e}$, where $\mathbf{e}$ denotes the identity matrix.

\begin{lem}\label{lem:tracefree}
Given any $\mathbf{u}_1,\mathbf{u}_2,\mathbf{u}_3\in\mathcal{M}_0$, one has
\begin{enumerate}
  \item[\rm(i)] $\mathbf{u}_1\mathbf{u}_2+\mathbf{u}_2\mathbf{u}_1={\rm tr}(\mathbf{u}_1\mathbf{u}_2)\cdot\mathbf{e}$;
  \item[\rm(ii)] $\mathbf{u}_1,\mathbf{u}_2,\mathbf{u}_3$ are linearly dependent if and only if ${\rm tr}(\mathbf{u}_1\mathbf{u}_2\mathbf{u}_3)=0$.
\end{enumerate}
\end{lem}

\begin{proof}
(i) Since $\mathbf{x}^\ast=-\mathbf{x}$ for $\mathbf{x}\in\mathcal{M}_0$, we have
\begin{align*}
{\rm tr}(\mathbf{u}_1\mathbf{u}_2)\cdot\mathbf{e}=\mathbf{u}_1\mathbf{u}_2+(\mathbf{u}_1\mathbf{u}_2)^\ast
=\mathbf{u}_1\mathbf{u}_2+(-\mathbf{u}_2)(-\mathbf{u}_1)=\mathbf{u}_1\mathbf{u}_2+\mathbf{u}_2\mathbf{u}_1.
\end{align*}

(ii) Just observe that the map
$(\mathcal{M}_0)^3\to\mathbb{C}$,
$(\mathbf{u},\mathbf{v},\mathbf{w})\mapsto{\rm tr}(\mathbf{u}\mathbf{v}\mathbf{w})$
is non-degenerate, trilinear and by (a), alternating.
\end{proof}

Given $\mathbf{x}_1,\ldots,\mathbf{x}_r\in\mathcal{M}$, let $\mathbb{C}\langle\mathbf{x}_1,\ldots,\mathbf{x}_r\rangle\subseteq\mathcal{M}$ denote the subspace they span; call the tuple $(\mathbf{x}_1,\ldots,\mathbf{x}_r)$ {\it reducible} if $\mathbf{x}_1,\ldots,\mathbf{x}_r$ have a common eigenvector, and {\it irreducible} otherwise.

Given $i_1,\ldots,i_r\in\{1,\ldots,4\}$, let $t_{i_1\cdots i_r}$ and $s_{i_1\cdots i_r}$ respectively denote the function $G^4\to\mathbb{C}$ sending $(\mathbf{x}_1,\ldots,\mathbf{x}_4)$ to ${\rm tr}(\mathbf{x}_{i_1}\cdots\mathbf{x}_{i_r})$ and ${\rm tr}(\check{\mathbf{x}}_{i_1}\cdots\check{\mathbf{x}}_{i_r})$, where
$\check{\mathbf{x}}_i=\mathbf{x}_i-(1/2){\rm tr}(\mathbf{x}_i)\cdot\mathbf{e}\in\mathcal{M}_0.$
Being invariant under simultaneous conjugation, they define functions on $\mathcal{X}(F_4)$.
Phrased alternatively, let
$F_4=\langle x_1,\ldots,x_4\mid-\rangle,$
then $t_{i_1\cdots i_r}$ sends
$\chi\in\mathcal{X}(F_4)$ to $\chi(x_{i_1}\cdots x_{i_r})$, and the $s_{i_1\cdots i_r}$'s for $r\ge 2$ can be expressed in terms of $t_{j_1\cdots j_k}$'s.

It is known that the coordinate ring $\mathbb{C}[\mathcal{X}(F_4)]$ is generated by the $t_{i_1\cdots i_r}$'s for $1\le i_1<\cdots<i_r\le 4$ and $1\le r\le 3$, and the {\it definitive relations} are
\begin{align}
{\rm(I)}: \ \ \ \  &2s_{i_1i_2i_3}s_{j_1j_2j_3}+\det\big((s_{i_aj_b})_{a,b=1}^3\big)=0,  \nonumber  \\
&\hspace{16mm}  1\le i_1<i_2<i_3\le 4,\ \ \ 1\le j_1<j_2<j_3\le 4;  \label{eq:typeI}  \\
{\rm(II)}: \ \ \ \   &s_{i1}s_{234}-s_{i2}s_{134}+s_{i3}s_{124}-s_{i4}s_{123}=0, \qquad  1\le i\le 4.  \label{eq:typeII}
\end{align}
One may refer to \cite[Theorem 3.1]{ABL18}.

For $t\in\mathbb{C}$, let $G(t)=\{\mathbf{x}\in{\rm SL}(2,\mathbb{C})\colon{\rm tr}(\mathbf{x})=t\}$.

For any $\mathbf{x},\mathbf{y}\in G(t)$, the following identities are true:
\begin{align}
\mathbf{x}^2&=t\cdot\mathbf{x}-\mathbf{e}, \qquad\qquad  \mathbf{x}^{-1}=t\mathbf{e}-\mathbf{x};  \label{eq:matrix-identity-2-1}  \\
\mathbf{x}\mathbf{y}\mathbf{x}&={\rm tr}(\mathbf{x}\mathbf{y})\cdot\mathbf{x}-\mathbf{y}^{-1}
={\rm tr}(\mathbf{x}\mathbf{y})\cdot\mathbf{x}+\mathbf{y}-t\mathbf{e}.  \label{eq:matrix-identity-2-2}
\end{align}
The first line is Cayley-Hamilton Theorem; the second line is deduced as
$$\mathbf{x}\mathbf{y}\mathbf{x}=(\mathbf{x}\mathbf{y})^2\mathbf{y}^{-1}
=({\rm tr}(\mathbf{x}\mathbf{y})\cdot\mathbf{x}\mathbf{y}-\mathbf{e})\mathbf{y}^{-1}
={\rm tr}(\mathbf{x}\mathbf{y})\cdot\mathbf{x}-\mathbf{y}^{-1}.$$

\begin{lem}\label{lem:irreducible}
For any $\mathbf{x},\mathbf{y}\in G(t)$, the following are equivalent to each other:
{\rm(i)} ${\rm tr}(\mathbf{x}\mathbf{y})\notin\{2,t^2-2\}$;
{\rm(ii)} $(\mathbf{x},\mathbf{y})$ is irreducible;
{\rm(iii)} $\mathbb{C}\langle\mathbf{e},\mathbf{x},\mathbf{y},\mathbf{x}\mathbf{y}\rangle=\mathcal{M}$.
\end{lem}

\begin{proof}
This is a part of \cite[Proposition 2.3.1]{Go09}.
For the equivalence between (i) and (ii), one may also refer to \cite[Lemma 2.2]{Ch22}.
\end{proof}

Let $s_0=t^2/2-2$, with $t$ fixed. On restriction to $G(t)^4$, we have
\begin{align}
s_{ii}=s_0,  \qquad  s_{ij}=t_{ij}-\frac{1}{2}t^2, \label{eq:sij-tij}  \\
s_{ijk}=t_{ijk}-\frac{t}{2}(t_{ij}+t_{jk}+t_{ik})+\frac{1}{2}t^3.   \label{eq:sijk-tijk}
\end{align}
An application of Lemma \ref{lem:tracefree} (i) yields
\begin{align}
\check{\mathbf{x}}_i^2=\frac{s_0}{2}\mathbf{e}, \qquad \check{\mathbf{x}}_i\check{\mathbf{x}}_j+\check{\mathbf{x}}_j\check{\mathbf{x}}_i=s_{ij}\mathbf{e}.  \label{eq:matrix-identity-3}
\end{align}

When simultaneously dealing with $4$ matrices, a probably frustrating problem is:
there are too many relations (20 in total) among the trace variables.
However, the situation is not so serious as it appears. In most cases, among the $16$ relations given by (\ref{eq:typeI}), only one needs to be handled.

\begin{lem}\label{lem:matrix}
{\rm(a)} Given $\mathsf{s}_{ij}$, $1\le i,j\le 3$ and $\mathsf{s}_{123}\ne 0$ such that $\mathsf{s}_{ii}=s_0$, $1\le i\le 3$ and $2\mathsf{s}_{123}^2+\det\big((\mathsf{s}_{ij})_{i,j=1}^3\big)=0$, up to simultaneous conjugacy there exists a unique $(\mathbf{x}_1,\mathbf{x}_2,\mathbf{x}_3)\in G(t)^3$ with $s_{ij}=\mathsf{s}_{ij}$ for all $i,j$ and $s_{123}=\mathsf{s}_{123}$.

{\rm(b)} Suppose $\mathbf{x}_1,\mathbf{x}_2,\mathbf{x}_3\in G(t)$, with $s_{123}=\mathsf{s}_{123}\ne 0$ and $s_{ij}=\mathsf{s}_{ij}$, $1\le i,j\le 3$.
Given $\mathsf{s}_{14},\ldots,\mathsf{s}_{44}$ and
$\mathsf{s}_{124}$, $\mathsf{s}_{134}$, $\mathsf{s}_{234}$ such that $\mathsf{s}_{44}=s_0$ and
$$\mathsf{s}_{i1}\mathsf{s}_{234}-\mathsf{s}_{i2}\mathsf{s}_{134}+\mathsf{s}_{i3}\mathsf{s}_{124}-\mathsf{s}_{i4}\mathsf{s}_{123}=0, \qquad
1\le i\le 4,$$
there exists a unique $\mathbf{x}_4\in G(t)$ with $s_{i4}=\mathsf{s}_{i4}$ and $s_{ij4}=\mathsf{s}_{ij4}$ for all $i,j$; actually,
$$\check{\mathbf{x}}_4=\frac{\mathsf{s}_{234}}{\mathsf{s}_{123}}\check{\mathbf{x}}_1-\frac{\mathsf{s}_{134}}{\mathsf{s}_{123}}\check{\mathbf{x}}_2
+\frac{\mathsf{s}_{124}}{\mathsf{s}_{123}}\check{\mathbf{x}}_3.$$
\end{lem}

\begin{proof}
(a) This was essentially proved in \cite[Section 5]{Go09}. Note that by Lemma \ref{lem:tracefree} (ii), $\check{\mathbf{x}}_1,\check{\mathbf{x}}_2,\check{\mathbf{x}}_3$ is linearly independent, so $(\mathbf{x}_1,\mathbf{x}_2,\mathbf{x}_3)$ is irreducible.

(b) Since $\mathsf{s}_{123}\ne 0$, so that $\mathbb{C}\langle\check{\mathbf{x}}_1,\check{\mathbf{x}}_2,\check{\mathbf{x}}_3\rangle=\mathcal{M}_0$,
there exists at most one $\mathbf{x}_4\in G(t)$ with prescribed ${\rm tr}(\check{\mathbf{x}}_i\check{\mathbf{x}}_4)$ for $1\le i\le3$.
Put $\check{\mathbf{x}}_4=\sum_{j=1}^3c_j\check{\mathbf{x}}_j$, with
$$c_1=\frac{\mathsf{s}_{234}}{\mathsf{s}_{123}}, \qquad c_2=-\frac{\mathsf{s}_{134}}{\mathsf{s}_{123}}, \qquad
c_3=\frac{\mathsf{s}_{124}}{\mathsf{s}_{123}}.$$
By Lemma \ref{lem:tracefree} (i), ${\rm tr}(\check{\mathbf{x}}_2\check{\mathbf{x}}_3\check{\mathbf{x}}_2)=0$, etc,
and ${\rm tr}(\check{\mathbf{x}}_1\check{\mathbf{x}}_3\check{\mathbf{x}}_2)=-{\rm tr}(\check{\mathbf{x}}_1\check{\mathbf{x}}_2\check{\mathbf{x}}_3)$.
Hence
\begin{align*}
s_{234}&={\rm tr}(\check{\mathbf{x}}_2\check{\mathbf{x}}_3\check{\mathbf{x}}_4)
={\sum}_{j=1}^3c_j{\rm tr}(\check{\mathbf{x}}_2\check{\mathbf{x}}_3\check{\mathbf{x}}_j)=\mathsf{s}_{234},   \\
s_{134}&={\rm tr}(\check{\mathbf{x}}_1\check{\mathbf{x}}_3\check{\mathbf{x}}_4)
={\sum}_{j=1}^3c_j{\rm tr}(\check{\mathbf{x}}_1\check{\mathbf{x}}_3\check{\mathbf{x}}_j)=\mathsf{s}_{134},   \\
s_{124}&={\rm tr}(\check{\mathbf{x}}_1\check{\mathbf{x}}_2\check{\mathbf{x}}_4)
={\sum}_{j=1}^3c_j{\rm tr}(\check{\mathbf{x}}_1\check{\mathbf{x}}_2\check{\mathbf{x}}_j)=\mathsf{s}_{124}.
\end{align*}
Moreover, for $i=1,2,3,4$,
\begin{align*}
s_{i4}&={\rm tr}(\check{\mathbf{x}}_i\check{\mathbf{x}}_4)={\sum}_{j=1}^3c_j{\rm tr}(\check{\mathbf{x}}_i\check{\mathbf{x}}_j)
=\frac{1}{\mathsf{s}_{123}}(\mathsf{s}_{i1}\mathsf{s}_{234}-\mathsf{s}_{i2}\mathsf{s}_{134}+\mathsf{s}_{i3}\mathsf{s}_{124})=\mathsf{s}_{i4}.
\end{align*}

Finally, $s_{44}=\mathsf{s}_{44}=s_0$ ensures $\mathbf{x}_4=\check{\mathbf{x}}_4+(t/2)\mathbf{e}\in G(t)$.
\end{proof}

The upshot is: any given $s_{i_1\cdots i_r}$'s can be realized by an irreducible quadruple $(\mathbf{x}_1,\mathbf{x}_2,\mathbf{x}_3,\mathbf{x}_4)\in G(t)^4$, which is unique up to conjugacy, as long as
$$2s_{123}^2=-\det\big((s_{ij})_{i,j=1}^3\big)\ne 0$$
and (\ref{eq:typeII}) holds.
The remaining $15$ relations in (\ref{eq:typeI}) are automatically guaranteed, as they are universal identities.

The following supplementary lemma will be useful:
\begin{lem}\label{lem:a1=a3}
Suppose $(\mathbf{x}_1,\mathbf{x}_2,\mathbf{x}_3)\in G(t)^3$ is irreducible,
$s_{123}=0$ and $s_{13}=\epsilon s_0$, with $\epsilon\in\{\pm1\}$. Then $\mathbf{x}_1\mathbf{x}_3=\mathbf{x}_3\mathbf{x}_1$.
Furthermore, $\mathbf{x}_3=\mathbf{x}_1^{\epsilon}$ if $t\ne\pm 2$.
\end{lem}

\begin{proof}
Assume $s_{13}=s_0$; for the case $s_{13}=-s_0$, just replace $\mathbf{x}_3$ by $\mathbf{x}_3^{-1}$.

The condition $s_{13}=s_0$ is equivalent to $t_{13}=t^2-2$.
By Lemma \ref{lem:irreducible}, $\mathbf{x}_1,\mathbf{x}_3$ have a common eigenvector, say $\xi$.
By Lemma \ref{lem:tracefree} (ii), $a_1\check{\mathbf{x}}_1+a_2\check{\mathbf{x}}_2+a_3\check{\mathbf{x}}_3=0$ for some $a_i$'s with at least one $a_i\ne0$.
If $a_2\ne0$, then $\mathbf{x}_2\in\mathbb{C}\langle\mathbf{e},\mathbf{x}_1,\mathbf{x}_3\rangle$, so $\xi$ is also an eigenvector of $\mathbf{x}_2$, contradicting the irreducibility of $(\mathbf{x}_1,\mathbf{x}_2,\mathbf{x}_3)$.
Hence $a_2=0$, implying $\mathbf{x}_1\mathbf{x}_3=\mathbf{x}_3\mathbf{x}_1$.

If $t\ne \pm2$, the up to conjugacy we may assume $\mathbf{x}_1=\mathbf{d}(\kappa)$ (the diagonal matrix with upper left entry $\kappa$) for some $\kappa$ with $\kappa+\kappa^{-1}=t$, then $\mathbf{x}_3$ must also be diagonal, say $\mathbf{x}_3=\mathbf{d}(\eta)$, with $\eta+\eta^{-1}=t$.
The condition $t_{13}=t^2-2$ forces $\eta=\kappa$. Hence $\mathbf{x}_3=\mathbf{x}_1$.
\end{proof}

\begin{figure}[h]
  \centering
  \includegraphics[width=4cm]{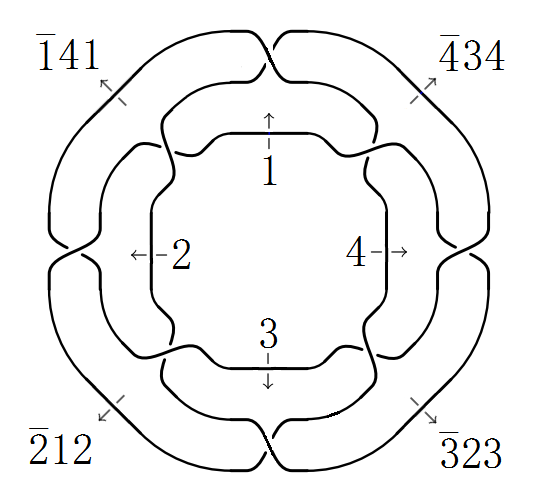}\\
  \caption{In the context of Wirtinger presentation, the arc labeled with $i$ gives rise to the generator $x_i$ for $\pi(T)$, and the one labeled with $\overline{1}41$ represents $x_1^{-1}x_4x_1\in\pi(T)$, and so on. This figure is modified from the one on
  \cite[Page 18]{PP20}.}\label{fig:rep}
\end{figure}

Referred to Figure \ref{fig:rep}, a Wirtinger presentation for $\pi(T)$ is given by
\begin{align}
\langle x_1,x_2,x_3,x_4\mid x_1x_2x_1^{-1}x_4x_1=x_2x_3x_2^{-1}x_1x_2=x_3x_4x_3^{-1}x_2x_3=x_4x_1x_4^{-1}x_3x_4\rangle.  \label{eq:presentation}
\end{align}

Let $f$ denote the rotation of angle $\pi/4$. Its induced automorphism on $\pi(T)$ sends $x_i$ to $x_{i+1}$.
The subscripts are understood as modulo $4$, as always done throughout.

The strategy of determining $\mathcal{X}^{\rm irr}(T)$ is outlined as follows.

Basically, $x_1,x_2,x_3,x_4$ belong to the same conjugacy class in $\pi(T)$.
Each character $\chi\in\mathcal{X}^{\rm irr}(T)$ is determined by
$$t=t_1,\  t_{12},\ t_{23},\ t_{34},\ t_{14},\ t_{13},\ t_{24}, \ t_{123},\ t_{124},\ t_{134},\ t_{234},$$
with $t_{i_1\cdots i_r}=\chi(x_{i_1}\cdots x_{i_r})$, so $\mathcal{X}^{\rm irr}(T)$ can be embedded into $\mathbb{C}^{11}$.

It is convenient to use $\vec{t}\in\mathbb{C}^{11}$ to denote an element of $\mathcal{X}^{\rm irr}(T)$.
The rotation $f$ induces an automorphism $f^\ast$ of $\mathcal{X}^{\rm irr}(T)$, by cyclically permuting the subscripts in the components of $\vec{t}$ except $t$.

Given $\mathbf{x}_1,\mathbf{x}_2,\mathbf{x}_3,\mathbf{x}_4\in G(t)$, there exists a (unique) representation $\rho:\pi(T)\to G$ sending $x_i$ to $\mathbf{x}_i$ if and only if $\mathbf{r}_1=\mathbf{r}_2=\mathbf{r}_3=\mathbf{r}_4$, where $$\mathbf{r}_i=\mathbf{x}_i\mathbf{x}_{i+1}\mathbf{x}_i^{-1}\mathbf{x}_{i-1}\mathbf{x}_i.$$
The key steps will be translating the matrix equations $\mathbf{r}_1=\mathbf{r}_2=\mathbf{r}_3=\mathbf{r}_4$ into {\it trace equations}, which are polynomial equations in terms of the $t_{i_1\cdots i_r}$'s, and finding all the solutions.

Let $S_{i_1i_2i_3}=(s_{i_ai_b})_{a,b=1}^3$. Let $S_\diamond=(s_{ij})_{i,j=1}^4$.

By Lemma \ref{lem:irreducible}, Lemma \ref{lem:matrix}, when $s_{123}\ne 0$,
we only need to treat
\begin{align}
2s_{123}^2+\det S_{123}=0   \label{eq:I}
\end{align}
and the type II relations given in (\ref{eq:typeII}), which are repeated here:
\begin{align}
s_{14}s_{123}+s_{12}s_{134}&=s_{11}s_{234}+s_{13}s_{124},  \label{eq:II-1}  \\
s_{24}s_{123}+s_{22}s_{134}&=s_{12}s_{234}+s_{23}s_{124},  \label{eq:II-2}  \\
s_{34}s_{123}+s_{23}s_{134}&=s_{13}s_{234}+s_{33}s_{124},  \label{eq:II-3}  \\
s_{44}s_{123}+s_{24}s_{134}&=s_{14}s_{234}+s_{34}s_{124}.  \label{eq:II-4}
\end{align}
The situation is similar when another $s_{ijk}\ne0$.

When all $s_{ijk}=0$, by looking into the trace equations carefully, we will always be able to find conditions to ensure $\mathbf{r}_1=\mathbf{r}_2=\mathbf{r}_3=\mathbf{r}_4$.

A particular consequence to (\ref{eq:II-1})--(\ref{eq:II-4}) is $\det S_\diamond=0$. Being a polynomial identity, it holds without the condition $s_{123}\ne 0$. One can refer to \cite[Theorem 3.1]{Do01} for this identity.
Usually it is helpful to use $\det S_\diamond=0$.

\section{From matrix equations to trace equations}

The irreducibility of $\rho$ is equivalent to that of $(\mathbf{x}_1,\mathbf{x}_2,\mathbf{x}_3,\mathbf{x}_4)$, which is always assumed.
In particular, $\mathbf{x}_i\ne\pm\mathbf{e}$ for all $i$. Actually, $(\mathbf{x}_1,\mathbf{x}_2,\mathbf{x}_3)$ is irreducible, since from (\ref{eq:presentation}) we can see that $x_1,x_2,x_3$ already generate $\pi(T)$. Similarly, each of $(\mathbf{x}_1,\mathbf{x}_2,\mathbf{x}_4)$, $(\mathbf{x}_1,\mathbf{x}_3,\mathbf{x}_4)$, $(\mathbf{x}_2,\mathbf{x}_3,\mathbf{x}_4)$ is irreducible.

\begin{lem}\label{lem:priori}
For any $i$, if $\mathbf{x}_i\mathbf{x}_{i+1}=\mathbf{x}_{i+1}\mathbf{x}_i$, then $\mathbf{x}_i=\mathbf{x}_{i+1}$, $\mathbf{x}_{i-1}=\mathbf{x}_{i+2}$, and ${\rm tr}(\mathbf{x}_i\mathbf{x}_{i-1})={\rm tr}(\mathbf{x}_{i+1}\mathbf{x}_{i+2})=1$.
\end{lem}

\begin{proof}
Without loss of generality, assume $\mathbf{x}_1\mathbf{x}_2=\mathbf{x}_2\mathbf{x}_1$.
Then
$\mathbf{r}_1=\mathbf{x}_2\mathbf{x}_4\mathbf{x}_1$, and $\mathbf{r}_2=\mathbf{x}_2\mathbf{x}_3\mathbf{x}_1$.
It follows from $\mathbf{r}_1=\mathbf{r}_2$ that $\mathbf{x}_4=\mathbf{x}_3$.

Consequently, $\mathbf{r}_3=\mathbf{x}_3\mathbf{x}_2\mathbf{x}_3$, and $\mathbf{r}_4=\mathbf{x}_3\mathbf{x}_1\mathbf{x}_3$.
From $\mathbf{r}_3=\mathbf{r}_4$ we see $\mathbf{x}_1=\mathbf{x}_2$.

Now $\mathbf{r}_2=\mathbf{x}_2\mathbf{x}_3\mathbf{x}_2$.
We can use (\ref{eq:matrix-identity-2-1}), (\ref{eq:matrix-identity-2-2}) to write $\mathbf{r}_2=\mathbf{r}_3$ as
$$0=\mathbf{r}_2-\mathbf{r}_3=(t_{23}\mathbf{x}_2-\mathbf{x}_3^{-1})-(t_{23}\mathbf{x}_3-\mathbf{x}_2^{-1})=(t_{23}-1)(\mathbf{x}_2-\mathbf{x}_3).$$
The irreducibility requires $\mathbf{x}_2\ne \mathbf{x}_3$, so $t_{23}=1$. And also $t_{14}=1$.
\end{proof}

Applying (\ref{eq:matrix-identity-2-2}) to $\mathbf{x}=\mathbf{x}_1$ and $\mathbf{y}=\mathbf{x}_2\mathbf{x}_1^{-1}\mathbf{x}_4$, we obtain
\begin{align*}
\mathbf{r}_1=\mathbf{x}_1\cdot\mathbf{x}_2\mathbf{x}_1^{-1}\mathbf{x}_4\cdot\mathbf{x}_1
={\rm tr}(\mathbf{x}_1\mathbf{x}_2\mathbf{x}_1^{-1}\mathbf{x}_4)\cdot\mathbf{x}_1-\mathbf{x}_4^{-1}\mathbf{x}_1\mathbf{x}_2^{-1}.
\end{align*}
Since
\begin{align*}
&{\rm tr}(\mathbf{x}_1\mathbf{x}_2\mathbf{x}_1^{-1}\mathbf{x}_4)
={\rm tr}(\mathbf{x}_1\mathbf{x}_2(t\mathbf{e}-\mathbf{x}_1)\mathbf{x}_4)  \\
=\ &t\cdot{\rm tr}(\mathbf{x}_1\mathbf{x}_2\mathbf{x}_4)-{\rm tr}((t_{12}\mathbf{x}_1+\mathbf{x}_2-t\mathbf{e})\mathbf{x}_4)
=tt_{124}-t_{12}t_{14}-t_{24}+t^2,
\end{align*}
we have
$$\mathbf{r}_1=(tt_{124}-t_{12}t_{14}-t_{24}+t^2)\mathbf{x}_1-\mathbf{x}_4^{-1}\mathbf{x}_1\mathbf{x}_2^{-1}.$$

Cyclically permuting the subscripts leads us to
\begin{align*}
\mathbf{r}_2&=(tt_{123}-t_{23}t_{12}-t_{13}+t^2)\mathbf{x}_2-\mathbf{x}_1^{-1}\mathbf{x}_2\mathbf{x}_3^{-1}, \\
\mathbf{r}_3&=(tt_{234}-t_{34}t_{23}-t_{24}+t^2)\mathbf{x}_3-\mathbf{x}_2^{-1}\mathbf{x}_3\mathbf{x}_4^{-1}, \\
\mathbf{r}_4&=(tt_{134}-t_{14}t_{34}-t_{13}+t^2)\mathbf{x}_4-\mathbf{x}_3^{-1}\mathbf{x}_4\mathbf{x}_1^{-1}.
\end{align*}

A necessary condition for $\mathbf{r}_1=\mathbf{r}_2$ is
\begin{align}
{\rm tr}(\mathbf{r}_1\mathbf{a})={\rm tr}(\mathbf{r}_2\mathbf{a}),  \qquad  \mathbf{a}\in\{\mathbf{e},\mathbf{x}_1,\mathbf{x}_2,\mathbf{x}_2\mathbf{x}_1\}.  \label{eq:necessary}
\end{align}

\begin{lem}
The condition {\rm(\ref{eq:necessary})} holds if and only if
\begin{align}
&(t^2-1)(t_{124}-t_{123})+t(t_{13}-t_{24})=t(t_{12}-1)(t_{14}-t_{23}),  \label{eq:tr0}   \\
&t(t_{12}+t^2-2)(t_{124}-t_{123})+(t_{12}+t^2-1)(t_{13}-t_{24})  \nonumber  \\
&\hspace{25mm} =(t_{12}^2+(t^2-1)t_{12}-t^2-1)(t_{14}-t_{23}),   \label{eq:tr1-tr2}  \\
&t(t^2-2-t_{12})(t_{124}+t_{123})+(t_{12}+1-t^2)(t_{13}+t_{24})  \nonumber  \\
&\hspace{25mm} =(1-t_{12}^2+(t^2-1)t_{12}-t^2)(t_{14}+t_{23}).   \label{eq:tr1+tr2}
\end{align}
\end{lem}

\begin{proof}
Since
\begin{align*}
&{\rm tr}(\mathbf{x}_4^{-1}\mathbf{x}_1\mathbf{x}_2^{-1})={\rm tr}((t\mathbf{e}-\mathbf{x}_4)\mathbf{x}_1(t\mathbf{e}-\mathbf{x}_2))  \\
=\ &{\rm tr}(t^2\mathbf{x}_1-t\mathbf{x}_1\mathbf{x}_2-t\mathbf{x}_4\mathbf{x}_1+\mathbf{x}_4\mathbf{x}_1\mathbf{x}_2)
=t^3-tt_{12}-tt_{14}+t_{124},
\end{align*}
we have
\begin{align*}
{\rm tr}(\mathbf{r}_1)&=(tt_{124}-t_{12}t_{14}-t_{24}+t^2)t-(t^3-tt_{12}-tt_{14}+t_{124})  \\
&=(t^2-1)t_{124}-tt_{12}t_{14}+t(t_{12}+t_{14}-t_{24}).
\end{align*}
Similarly (or by permuting the subscripts),
$${\rm tr}(\mathbf{r}_2)=(t^2-1)t_{123}-tt_{23}t_{12}+t(t_{23}+t_{12}-t_{13}).$$
Hence ${\rm tr}(\mathbf{r}_1)={\rm tr}(\mathbf{r}_2)$ is equivalent to (\ref{eq:tr0}).

\medskip

Let $t_{i\overline{k}}={\rm tr}(\mathbf{x}_i\mathbf{x}_k^{-1})={\rm tr}(\mathbf{x}_i(t\mathbf{e}-\mathbf{x}_k))=t^2-t_{ik}$.
Then
\begin{align*}
&{\rm tr}(\mathbf{x}_4^{-1}\mathbf{x}_1\mathbf{x}_2^{-1}\mathbf{x}_1)
={\rm tr}(\mathbf{x}_4^{-1}(t_{1\overline{2}}\mathbf{x}_1-\mathbf{x}_2))=t_{1\overline{2}}t_{1\overline{4}}-t_{2\overline{4}}  \\
=\ &(t^2-t_{12})(t^2-t_{14})-(t^2-t_{24})=t^4-t^2(t_{12}+t_{14}+1)+t_{12}t_{14}+t_{24},  \\
&{\rm tr}(\mathbf{x}_1^{-1}\mathbf{x}_2\mathbf{x}_3^{-1}\mathbf{x}_2)
={\rm tr}(\mathbf{x}_1^{-1}(t_{2\overline{3}}\mathbf{x}_2-\mathbf{x}_3))=t_{2\overline{3}}t_{2\overline{1}}-t_{3\overline{1}}  \\
=\ &(t^2-t_{23})(t^2-t_{12})-(t^2-t_{13})=t^4-t^2(t_{12}+t_{23}+1)+t_{12}t_{23}+t_{13}.
\end{align*}

We can compute
\begin{align*}
{\rm tr}(\mathbf{r}_1\mathbf{x}_1)-{\rm tr}(\mathbf{r}_2\mathbf{x}_1)
&=(tt_{124}-t_{12}t_{14}-t_{24}+t^2)(t^2-2)-{\rm tr}(\mathbf{x}_4^{-1}\mathbf{x}_1\mathbf{x}_2^{-1}\mathbf{x}_1)  \\
&\ \ \ \ -(tt_{123}-t_{23}t_{12}-t_{13}+t^2)t_{12}+t_{2\overline{3}}  \\
&=(t^2-2)tt_{124}-tt_{12}t_{123}+t_{12}t_{13}+(1-t^2)(t_{12}t_{14}+t_{24}) \\
&\ \ \ +(t_{12}^2-1)t_{23}+t^2t_{14},   \\
{\rm tr}(\mathbf{r}_1\mathbf{x}_2)-{\rm tr}(\mathbf{r}_2\mathbf{x}_2)
&=(tt_{123}-t_{12}t_{23}-t_{13}+t^2)t_{12}-t_{2\overline{3}} \\
&\ \ \ \ -(tt_{123}-t_{23}t_{12}-t_{13}+t^2)(t^2-2)+{\rm tr}(\mathbf{x}_1^{-1}\mathbf{x}_2\mathbf{x}_3^{-1}\mathbf{x}_2)  \\
&=(t^2-2)tt_{123}-tt_{12}t_{124}+t_{12}t_{24}+(1-t^2)(t_{12}t_{23}+t_{13}) \\
&\ \ \ +(t_{12}^2-1)t_{14}+t^2t_{23}.
\end{align*}
Hence ${\rm tr}(\mathbf{r}_1\mathbf{x}_1)={\rm tr}(\mathbf{r}_2\mathbf{x}_1)$ is equivalent to
\begin{align}
&(t^2-2)tt_{124}-tt_{12}t_{123}+t_{12}t_{13}+(1-t^2)t_{24}  \nonumber  \\
=\ &(1-t_{12}^2)t_{23}+(t^2-1)t_{12}t_{14}-t^2t_{14}, \label{eq:tr1}
\end{align}
and ${\rm tr}(\mathbf{r}_1\mathbf{x}_2)={\rm tr}(\mathbf{r}_2\mathbf{x}_2)$ is equivalent to
\begin{align}
&(t^2-2)tt_{123}-tt_{12}t_{124}+t_{12}t_{24}+(1-t^2)t_{13}  \nonumber  \\
=\ &(1-t_{12}^2)t_{14}+(t^2-1)t_{12}t_{23}-t^2t_{23}. \label{eq:tr2}
\end{align}

Furthermore, using
\begin{align*}
{\rm tr}(\mathbf{x}_i^2\mathbf{x}_k)&={\rm tr}((t\mathbf{x}_i-\mathbf{e})\mathbf{x}_k)=tt_{ik}-t,   \\
{\rm tr}(\mathbf{x}_i^2\mathbf{x}_k^{-1})&={\rm tr}((t\mathbf{x}_i-\mathbf{e})\mathbf{x}_k^{-1})=tt_{i\overline{k}}-t=t(t^2-t_{ik})-t,
\end{align*}
we compute
\begin{align*}
&{\rm tr}(\mathbf{r}_1\mathbf{x}_2\mathbf{x}_1)-{\rm tr}(\mathbf{r}_2\mathbf{x}_2\mathbf{x}_1)  \\
=\ &(tt_{124}-t_{12}t_{14}-t_{24}+t^2){\rm tr}(\mathbf{x}_1^2\mathbf{x}_2)-{\rm tr}(\mathbf{x}_1^2\mathbf{x}_4^{-1})  \\
&-(tt_{123}-t_{12}t_{23}-t_{13}+t^2){\rm tr}(\mathbf{x}_2^2\mathbf{x}_1)+{\rm tr}(\mathbf{x}_2^2\mathbf{x}_3^{-1})  \\
=\ &t^2(t_{12}-1)(t_{124}-t_{123})+t(t_{12}-1)(t_{13}-t_{24})+t(1+t_{12}-t_{12}^2)(t_{14}-t_{23}).
\end{align*}
Hence
${\rm tr}(\mathbf{r}_1\mathbf{x}_2\mathbf{x}_1)={\rm tr}(\mathbf{r}_2\mathbf{x}_2\mathbf{x}_1)$ if and only if
\begin{align}
t^2(t_{12}-1)(t_{124}-t_{123})+t(t_{12}-1)(t_{13}-t_{24})=t(t_{12}^2-t_{12}-1)(t_{14}-t_{23}).  \label{eq:tr3}
\end{align}

Taking the difference between (\ref{eq:tr1}) and (\ref{eq:tr2}) yields (\ref{eq:tr1-tr2});
summing (\ref{eq:tr1}) and (\ref{eq:tr2}) yields (\ref{eq:tr1+tr2}).
Thus, (\ref{eq:tr1}), (\ref{eq:tr2}) can be replaced by (\ref{eq:tr1-tr2}), (\ref{eq:tr1+tr2}).

Finally, observe that (\ref{eq:tr3}) is redundant: it follows from (\ref{eq:tr0}) and (\ref{eq:tr1-tr2}).
\end{proof}

\begin{rmk}\label{rmk:trace-equation}
\rm (i) If $t_{12}\ne 1$, then from (\ref{eq:tr0}), (\ref{eq:tr1-tr2}) we solve
\begin{align}
t_{124}-t_{123}&=\Big(t+\frac{t}{1-t_{12}}\Big)(t_{14}-t_{23}), \label{eq:trace-eq-1}  \\
t_{13}-t_{24}&=\Big(t_{12}-t^2+\frac{t^2-1}{t_{12}-1}\Big)(t_{14}-t_{23}).  \label{eq:trace-eq-2}
\end{align}
So (\ref{eq:necessary}) is alternatively equivalent to (\ref{eq:tr1+tr2}), (\ref{eq:trace-eq-1}), (\ref{eq:trace-eq-2}).

(ii) If $t_{12}=1$, then (\ref{eq:tr0}) becomes
\begin{align}
(t^2-1)(t_{124}-t_{123})+t(t_{13}-t_{24})=0,  \label{eq:trace-eq-1'}
\end{align}
and (\ref{eq:tr1-tr2}) can be replaced by $t_{14}=t_{23}$; summing $t$ times of (\ref{eq:trace-eq-1'}) with (\ref{eq:tr1+tr2}) and dividing by $2$, we obtain
\begin{align*}
(t^2-2)tt_{124}-tt_{123}+t_{13}+(1-t^2)t_{24}+t_{14}=0.
\end{align*}

Sometimes it is useful to remember that (\ref{eq:tr1+tr2}) still holds when $t_{12}=1$.

(iii) If $t_{12}\ne 2,t^2-2$, then by Lemma \ref{lem:irreducible},
$\mathbb{C}\langle\mathbf{e},\mathbf{x}_1,\mathbf{x}_2,\mathbf{x}_2\mathbf{x}_1\rangle=\mathcal{M}$,
so (\ref{eq:necessary}) is also sufficient for $\mathbf{r}_1=\mathbf{r}_2$.
When $t_{12}\in\{2,t^2-2\}$, we will seek other conditions for $\mathbf{r}_1=\mathbf{r}_2$.

(iv) Bear in mind that everything is rotationally symmetric.
\end{rmk}

There turn out to be various cases to be discussed separately.
In each case, we will solve the trace equations (i.e., (\ref{eq:tr0})--(\ref{eq:tr1+tr2}) and those obtained by cyclically permuting the subscripts, as well as the equations equivalent to them), and verify the definitive relations if at least one of
$s_{123}$, $s_{124}$, $s_{134}$, $s_{234}$ is nonzero.

\section{Solving the trace equations}

Throughout, let $\alpha=t^2-1$. Let $x_{i,i+1}=t_{i,i+1}-1$.
Keep (\ref{eq:sij-tij}), (\ref{eq:sijk-tijk}) in mind; in particular, $s_{ij}-s_0=t_{ij}+2-t^2$, and $s_{ij}+s_0=t_{ij}-2$.
We shall frequently switch between $\alpha$, $x_{i,i+1}$, $s_{ij}$ and $t^2$, $t_{i,i+1}$, $t_{ij}$, respectively.

\subsection{$x_{12}x_{23}x_{34}x_{14}\ne0$}

List (\ref{eq:trace-eq-1}) and the equations obtained by rotations as
\begin{align}
t_{124}-t_{123}&=t\Big(1-\frac{1}{x_{12}}\Big)(x_{14}-x_{23}),  \label{eq:1-1-1}  \\
t_{123}-t_{234}&=t\Big(1-\frac{1}{x_{23}}\Big)(x_{12}-x_{34}),  \label{eq:1-1-2}  \\
t_{234}-t_{134}&=t\Big(1-\frac{1}{x_{34}}\Big)(x_{23}-x_{14}),  \label{eq:1-1-3}  \\
t_{134}-t_{124}&=t\Big(1-\frac{1}{x_{14}}\Big)(x_{34}-x_{12}).   \label{eq:1-1-4}
\end{align}
List (\ref{eq:trace-eq-2}) and the equations obtained by rotations as
\begin{align}
t_{13}-t_{24}&=\Big(x_{12}-\alpha+\frac{\alpha}{x_{12}}\Big)(x_{14}-x_{23}),  \label{eq:1-2-1}  \\
t_{24}-t_{13}&=\Big(x_{23}-\alpha+\frac{\alpha}{x_{23}}\Big)(x_{12}-x_{34}),  \label{eq:1-2-2}  \\
t_{13}-t_{24}&=\Big(x_{34}-\alpha+\frac{\alpha}{x_{34}}\Big)(x_{23}-x_{14}),  \label{eq:1-2-3}  \\
t_{24}-t_{13}&=\Big(x_{14}-\alpha+\frac{\alpha}{x_{14}}\Big)(x_{34}-x_{12}).  \label{eq:1-2-4}
\end{align}
List (\ref{eq:tr1+tr2}) and the equations obtained by rotations as
\begin{align}
t(t_{12}+1-\alpha)(t_{124}+t_{123})+(\alpha-t_{12})(t_{13}+t_{24})
=(t_{12}^2-\alpha t_{12}+\alpha)(t_{14}+t_{23}),  \label{eq:1-3-1} \\
t(t_{23}+1-\alpha)(t_{123}+t_{234})+(\alpha-t_{23})(t_{13}+t_{24})
=(t_{23}^2-\alpha t_{23}+\alpha)(t_{12}+t_{34}),  \label{eq:1-3-2}  \\
t(t_{34}+1-\alpha)(t_{234}+t_{134})+(\alpha-t_{34})(t_{13}+t_{24})
=(t_{34}^2-\alpha t_{34}+\alpha)(t_{14}+t_{23}),  \label{eq:1-3-3}  \\
t(t_{14}+1-\alpha)(t_{134}+t_{124})+(\alpha-t_{14})(t_{13}+t_{24})
=(t_{14}^2-\alpha t_{14}+\alpha)(t_{12}+t_{34}).  \label{eq:1-3-4}
\end{align}

\begin{ass}
If $t_{13}=t_{24}$, then $(x_{12}-x_{34})(x_{14}-x_{23})=0$.
\end{ass}

\begin{proof}
Assume $t_{13}=t_{24}$ but $(x_{12}-x_{34})(x_{14}-x_{23})\ne0$.

From (\ref{eq:1-2-1})--(\ref{eq:1-2-4}) we see that $x_{12},x_{34}$ are the two distinct roots of $x^2-\alpha x+\alpha=0$, so are $x_{14},x_{23}$.
This in particular forces $\alpha^2\ne 4\alpha$.

By symmetry we may just assume $x_{14}=x_{12}$ and $x_{34}=x_{23}$.

Comparing (\ref{eq:1-1-2}) with (\ref{eq:1-1-3}), we obtain $t_{134}=t_{123}$, so $s_{134}=s_{123}$.
Summing (\ref{eq:1-1-1}) and (\ref{eq:1-1-2}) yields
$t_{124}-t_{234}=t(x_{12}-x_{23})$, which implies $s_{124}=s_{234}$.

The equations (\ref{eq:1-3-1}), (\ref{eq:1-3-2}) respectively become
\begin{align*}
(t_{12}+1-\alpha)t(t_{124}+t_{123})&=(t_{12}^2-\alpha t_{12}+\alpha)(t_{14}+t_{23})+2(t_{12}-\alpha)t_{13},   \\
(t_{23}+1-\alpha)t(t_{123}+t_{234})&=(t_{23}^2-\alpha t_{23}+\alpha)(t_{12}+t_{34})+2(t_{23}-\alpha)t_{13}.
\end{align*}
Subtracting $(t_{12}+1-\alpha)$ times of the second row from $(t_{23}+1-\alpha)$ times of the first row yields
$$(\alpha^2-\alpha-6-2t_{13})(x_{12}-x_{23})=(\alpha-4)t(t_{124}-t_{234})=(\alpha-4)t^2(x_{12}-x_{23}).$$
Hence $t_{13}=\alpha-1$, implying $s_{13}=(\alpha-3)/2=s_0$.

Taking the difference between (\ref{eq:II-1}) and (\ref{eq:II-3}), we obtain
$$2(x_{12}-x_{23})s_{123}=(s_{11}-s_{13})(s_{234}-s_{124})=0.$$
So $s_{123}=0$, which by (\ref{eq:I}) implies $\det S_{123}=0$. Then we can deduce $\alpha=3$, i.e. $t\in\{\pm2\}$. Now $t_{13}=2$, so by Lemma \ref{lem:a1=a3}, $\mathbf{x}_1\mathbf{x}_3=\mathbf{x}_3\mathbf{x}_1$. Moreover,
$$0=(s_{24}+s_{22})s_{123}\stackrel{(\ref{eq:II-2})}=(s_{12}+s_{23})s_{234}=(x_{12}+x_{23}-2)s_{234}=s_{234}.$$

Since $x_{23}^2-3x_{23}+3=0$, one has $x_{23}\ne 1$, so that $t_{23}\ne 2=t^2-2$. By Lemma \ref{lem:irreducible},
$\mathbb{C}\langle\mathbf{e},\mathbf{x}_2,\mathbf{x}_3,\mathbf{x}_2\mathbf{x}_3\rangle=\mathcal{M}$.
As we have seen, ${\rm tr}(\mathbf{x}_2\mathbf{a})={\rm tr}(\mathbf{x}_4\mathbf{a})$ for
$\mathbf{a}\in\{\mathbf{e},\mathbf{x}_2,\mathbf{x}_3,\mathbf{x}_2\mathbf{x}_3\}$, so $\mathbf{x}_2=\mathbf{x}_4$.

Then $\mathbf{r}_2=\mathbf{r}_4$ implies $\mathbf{x}_3\mathbf{x}_2^{-1}\mathbf{x}_1=\mathbf{x}_1\mathbf{x}_2^{-1}\mathbf{x}_3$, so that
$\mathbf{x}_2$ commutes with $\mathbf{x}_1^{-1}\mathbf{x}_3$. This contradicts the irreducibility of $(\mathbf{x}_1,\mathbf{x}_2,\mathbf{x}_3)$.
\end{proof}

Therefore, in the case-by-case discussion below, we do not consider the case $t_{13}=t_{24}$, $(x_{12}-x_{34})(x_{14}-x_{23})\ne0$.

\subsubsection{$t_{13}\ne t_{24}$}

By (\ref{eq:1-2-1}), $t_{14}\ne t_{23}$; by (\ref{eq:1-2-2}), $t_{12}\ne t_{34}$.

Summing (\ref{eq:1-2-1}) and (\ref{eq:1-2-3}) yields
\begin{align}
(x_{12}+x_{34})\Big(1+\frac{\alpha}{x_{12}x_{34}}\Big)=2\alpha;   \label{eq:12-34}
\end{align}
summing (\ref{eq:1-2-2}) and (\ref{eq:1-2-4}) yields
\begin{align}
(x_{14}+x_{23})\Big(1+\frac{\alpha}{x_{14}x_{23}}\Big)=2\alpha.   \label{eq:14-23}
\end{align}

\begin{ass}
$\{x_{12},x_{34}\}=\{x_{14},x_{23}\}.$
\end{ass}

\begin{proof}
When $t\ne 0$, summing (\ref{eq:1-1-1})--(\ref{eq:1-1-4}) together, we obtain
\begin{align*}
0&=\Big(\frac{1}{x_{34}}-\frac{1}{x_{12}}\Big)(x_{14}-x_{23})+\Big(\frac{1}{x_{14}}-\frac{1}{x_{23}}\Big)(x_{12}-x_{34})  \nonumber  \\
&=\frac{1}{x_{12}x_{34}x_{14}x_{23}}(x_{14}x_{23}-x_{12}x_{34})(x_{12}-x_{34})(x_{14}-x_{23}),
\end{align*}
implying $x_{14}x_{23}=x_{12}x_{34}$. Then by (\ref{eq:12-34}), (\ref{eq:14-23}), $x_{12}+x_{34}=x_{14}+x_{23}$.
Hence $\{x_{12},x_{34}\}=\{x_{14},x_{23}\}$.

When $t=0$, subtracting (\ref{eq:1-3-3}) from (\ref{eq:1-3-1}) and dividing by $t_{34}-t_{12}$, we obtain
$$t_{13}+t_{24}+(t_{12}+t_{34})(t_{14}+t_{23})+t_{14}+t_{23}=0;$$
subtracting (\ref{eq:1-3-4}) from (\ref{eq:1-3-2}) and dividing by $t_{14}-t_{23}$, we obtain
$$t_{13}+t_{24}+(t_{14}+t_{23})(t_{12}+t_{34})+t_{12}+t_{34}=0.$$
Hence $x_{12}+x_{34}=x_{14}+x_{23}$, and then by (\ref{eq:12-34}), (\ref{eq:14-23}), $x_{12}x_{34}=x_{14}x_{23}$.
Thus, also $\{x_{12},x_{34}\}=\{x_{14},x_{23}\}$.
\end{proof}

Let $\mathcal{X}_{1,1}$ (resp. $\mathcal{X}_{1,2}$) denote the subset of $\mathcal{X}^{\rm irr}(T)$ consisting of characters with
$x_{12}x_{23}x_{34}x_{14}\ne 0$, $t_{13}\ne t_{24}$ as well as $x_{14}=x_{12}$, $x_{34}=x_{23}$ (resp. as well as $x_{23}=x_{12}$, $x_{34}=x_{14}$).

\begin{rmk}
\rm Due to that $f^\ast:\mathcal{X}^{\rm irr}(T)\to\mathcal{X}^{\rm irr}(T)$ interchanges $\mathcal{X}_{1,1}$ with $\mathcal{X}_{1,2}$,
we only deal with $\mathcal{X}_{1,1}$, replacing $x_{14},x_{34}$ respectively by $x_{12},x_{23}$ everywhere.

At the end of this subsection, we will directly state the result for $\mathcal{X}_{1,2}$.
\end{rmk}

Let
$\mu=x_{12}+x_{23}$, $\nu=x_{12}x_{23}.$
Rewrite (\ref{eq:12-34}) as
\begin{align}
(2\alpha-\mu)\nu=\alpha\mu, \qquad  \text{or\ \ \ \ \ } \mu\nu=\alpha(2\nu-\mu).    \label{eq:mu-nu}
\end{align}

Half of the sum of (\ref{eq:1-2-1}) and (\ref{eq:1-2-3}) leads to
\begin{align}
t_{13}-t_{24}&=\frac{1}{2}\big(1-\frac{\alpha}{\nu}\big)(x_{12}-x_{23})^2=\frac{1}{2}\big(1-\frac{\alpha}{\nu}\big)(\mu^2-4\nu) \nonumber  \\
&=\frac{1}{2}(\mu^2-4\nu-\mu(2\alpha-\mu)+4\alpha)=\mu^2-2\nu+\alpha(2-\mu),  \label{eq:t13-minus-t24}
\end{align}
where (\ref{eq:mu-nu}) is used in the first equality in the second line.

The sum of (\ref{eq:1-1-2}) and (\ref{eq:1-1-3}) yields
$$t_{134}=t_{123},$$
and the sum of (\ref{eq:1-1-1}) and (\ref{eq:1-1-2}) yields
\begin{align}
t_{124}-t_{234}=t\Big(2-\frac{1}{x_{12}}-\frac{1}{x_{23}}\Big)(x_{12}-x_{23}).  \label{eq:t124-t234}
\end{align}

Let $\gamma_1=t_{12}+1-\alpha$, $\gamma_2=t_{23}+1-\alpha$. Then
$$\gamma_1\gamma_2=\nu+(2-\alpha)\mu+(2-\alpha)^2.$$
Rewrite (\ref{eq:1-3-1}), (\ref{eq:1-3-2}) as
\begin{align}
\gamma_1t(t_{124}+t_{123})&=(\gamma_1^2+(\alpha-2)\gamma_1+1)(\mu+2)+(\gamma_1-1)(t_{13}+t_{24}),  \label{eq:last-1}  \\
\gamma_2t(t_{123}+t_{234})&=(\gamma_2^2+(\alpha-2)\gamma_2+1)(\mu+2)+(\gamma_2-1)(t_{13}+t_{24}).   \label{eq:last-2}
\end{align}
Subtracting $\gamma_1$ times of (\ref{eq:last-2}) from $\gamma_2$ times of (\ref{eq:last-1}) yields
\begin{align}
\gamma_1\gamma_2t(t_{124}-t_{234})=(x_{12}-x_{23})\big((\gamma_1\gamma_2-1)(\mu+2)+t_{13}+t_{24}\big).    \label{eq:subtract-1}
\end{align}

\begin{ass}
$\alpha\ne0$.
\end{ass}

\begin{proof}
Assume $\alpha=0$, i.e. $t^2=1$. Then by (\ref{eq:mu-nu}), $\mu=0$.
Let $x=x_{12}=x_{14}$. Then
$x_{34}=x_{23}=-x$; by (\ref{eq:t13-minus-t24}), (\ref{eq:t124-t234}),
$$t_{13}-t_{24}=2x^2,  \qquad  t_{124}-t_{234}=4tx.$$
Hence $\gamma_1\gamma_2=\nu+4=4-x^2$, and (\ref{eq:subtract-1}) becomes
\begin{align*}
0&=(4-x^2)t(t_{124}-t_{234})-2x(6-2x^2+t_{13}+t_{24})   \\
&=(4-x^2)\cdot 4t^2x-2x(6-4x^2+2t_{13})=4x+4x^3-4xt_{13},
\end{align*}
so $t_{13}=1+x^2$.
However, direct computation shows
$2s_{123}^2+\det S_{123}=6x^2\ne 0$, contradicting (\ref{eq:I}).
\end{proof}

Now (\ref{eq:mu-nu}) is equivalent to
$$\frac{\mu}{\nu}=2-\frac{\mu}{\alpha}.$$

The difference of (\ref{eq:1-1-1}) and (\ref{eq:1-1-2}) gives rise to
\begin{align}
t_{124}+t_{234}=2t_{123}+t\Big(\frac{\mu^2}{\nu}-4\Big)=2t_{123}+t\Big(2\mu-\frac{\mu^2}{\alpha}-4\Big).  \label{eq:t124+t234}
\end{align}

Rewrite (\ref{eq:t124-t234}) as
\begin{align}
t_{124}-t_{234}=t\big(2-\frac{\mu}{\nu}\big)(x_{12}-x_{23})=\frac{t\mu}{\alpha}(x_{12}-x_{23}),   \label{eq:t124-minus-t234}
\end{align}
which combined with (\ref{eq:subtract-1}) and (\ref{eq:mu-nu}) imply
\begin{align*}
t_{13}+t_{24}&=\gamma_1\gamma_2\big(1+\frac{1}{\alpha}\big)\mu-\big(\gamma_1\gamma_2-1\big)(\mu+2)   \\
&=\big(\frac{2}{\alpha}-1\big)\mu^2+\big(3\alpha-8+\frac{4}{\alpha}\big)\mu-2(\alpha^2-4\alpha+3).
\end{align*}
From this and (\ref{eq:t13-minus-t24}) we can solve
\begin{align}
t_{13}&=\frac{1}{\alpha}\mu^2-\nu+\big(\alpha-4+\frac{2}{\alpha}\big)\mu-\alpha^2+5\alpha-3,  \label{eq:solve-t13-1}  \\
t_{24}&=\big(\frac{1}{\alpha}-1\big)\mu^2+\nu+\big(2\alpha-4+\frac{2}{\alpha}\big)\mu-\alpha^2+3\alpha-3.  \label{eq:solve-t24-1}
\end{align}

We can replace (\ref{eq:II-1}) and (\ref{eq:II-3}) by their difference and sum, namely,
\begin{align}
2(s_{12}-s_{23})s_{123}&=(s_{11}-s_{13})(s_{234}-s_{124}), \label{eq:typeII-1-3} \\
2(s_{12}+s_{23})s_{123}&=(s_{11}+s_{13})(s_{234}+s_{124}). \label{eq:typeII-1+3}
\end{align}
Then (\ref{eq:II-2}), (\ref{eq:II-4}) are both equivalent to
\begin{align}
(t_{24}-2)s_{123}=s_{12}s_{234}+s_{23}s_{124};   \label{eq:typeII-2}
\end{align}
when $s_{123}\ne 0$, it is in turn equivalent to
\begin{align}
t_{24}-2=\frac{(s_{12}+s_{23})^2}{s_{11}+s_{13}}+\frac{(s_{12}-s_{23})^2}{s_{11}-s_{13}}.   \label{eq:II-2'}
\end{align}

From (\ref{eq:typeII-1-3}) we obtain
\begin{align*}
s_{123}=\frac{(s_{11}-s_{13})(t_{234}-t_{124}+t(x_{12}-x_{23}))}{2(x_{12}-x_{23})}=\frac{t}{2}(\alpha-1-t_{13})\big(1-\frac{\mu}{\alpha}\big).
\end{align*}
Hence by (\ref{eq:sijk-tijk}),
\begin{align}
t_{123}&=s_{123}+\frac{t}{2}(t_{13}+t_{12}+t_{23}-t^2)  \nonumber  \\
&=\frac{t}{2}(\alpha-1-t_{13})(1-\frac{\mu}{\alpha})+\frac{t}{2}(t_{13}+\mu+1-\alpha) =\frac{t\mu(t_{13}+1)}{2\alpha}.  \label{eq:1-t123}
\end{align}

Subtracting $(\gamma_1-1)$ times of (\ref{eq:last-2}) from $(\gamma_2-1)$ times of (\ref{eq:last-1}), we obtain
\begin{align*}
&\big(\nu+(1-\alpha)\mu+(\alpha-1)(\alpha-2)-1\big)(\mu+2)  \\
=\ &(\nu+(1-\alpha)\mu+(\alpha-1)(\alpha-2))\frac{t(t_{124}-t_{234})}{x_{12}-x_{23}}-tt_{123}+t\frac{x_{23}t_{124}-x_{12}t_{234}}{x_{12}-x_{23}}  \\
=\ &(\nu+(1-\alpha)\mu+(\alpha-1)(\alpha-2))\frac{t^2\mu}{\alpha}-\frac{t^2\mu(t_{13}+1)}{\alpha}+t^2\big(\frac{\mu^2}{\alpha}+2-\mu\big), \end{align*}
where in the last line we have applied (\ref{eq:t124+t234}) and (\ref{eq:t124-minus-t234}) to compute
\begin{align*}
\frac{x_{23}t_{124}-x_{12}t_{234}}{x_{12}-x_{23}}
&=-\frac{t_{124}+t_{234}}{2}+\frac{x_{12}+x_{23}}{2}\cdot\frac{t_{124}-t_{234}}{x_{12}-x_{23}}  \\
&=-t_{123}-t\big(\mu-\frac{\mu^2}{2\alpha}-2\big)+\frac{t\mu^2}{2\alpha}=-t_{123}+t\big(\frac{\mu^2}{\alpha}+2-\mu\big).
\end{align*}
Consequently,
\begin{align*}
&2\big(\mu^2+(\alpha^2-4\alpha)\mu+4\alpha^2-\alpha^3\big)=(\alpha+1)\mu(t_{13}-1) \\
=\ &(\alpha+1)\mu\Big(\frac{1}{\alpha}\mu^2-\nu+\big(\alpha-4+\frac{2}{\alpha}\big)\mu-\alpha^2+5\alpha-4\Big)  \\
=\ &\frac{(\alpha+1)\mu}{\alpha(\mu-2\alpha)}(\mu-2\alpha+2)(\mu^2+(\alpha^2-4\alpha)\mu+4\alpha^2-\alpha^3\big);
\end{align*}
in the second line, $t_{13}$ is substituted via (\ref{eq:solve-t13-1}), and in the third line, $\nu$ is substituted via (\ref{eq:mu-nu}).
Hence
\begin{align*}
\big((\alpha+1)\mu^2-2(\alpha^2+\alpha-1)\mu+4\alpha^2\big)\big(\mu^2+(\alpha^2-4\alpha)\mu+4\alpha^2-\alpha^3\big)=0.
\end{align*}

\begin{ass}
Each element of $\mathcal{X}_{1,1}$ satisfies
\begin{align}
\mu^2+(\alpha^2-4\alpha)\mu+4\alpha^2-\alpha^3=0.   \label{eq:component-1}
\end{align}
\end{ass}

\begin{proof}
Assume there exists $\chi_0\in\mathcal{X}_{1,1}$ not satisfying (\ref{eq:component-1}).
Then for $\chi_0$,
\begin{align*}
(\alpha+1)\mu^2-2(\alpha^2+\alpha-1)\mu+4\alpha^2=0.
\end{align*}
Note that $\alpha\ne-1$: otherwise $\mu=-2$, contradicting (\ref{eq:mu-nu}).
Hence
$$g(\mu):=\mu^2-\frac{2(\alpha^2+\alpha-1)}{\alpha+1}\mu+\frac{4\alpha^2}{\alpha+1}=0$$
near $\chi_0$.
We can rewrite $\nu$, $t_{13}$, $t_{24}$, $s_{12}+s_{23}$, $s_{12}s_{23}$ as elements of $\mathbb{Q}(\alpha)[\mu]$,
and when $s_{123}\ne 0$, reformulate (\ref{eq:II-2'}) as $h(\mu)=0$ for some $h(\mu)\in\mathbb{Q}(\alpha)[\mu]$.

We could disprove (\ref{eq:II-2'}) by direct computation, but it turns out to be very tedious. Here is a better approach.
Let $\mathcal{C}$ be an irreducible component of $\mathcal{X}^{\rm irr}(T)$ containing $\chi_0$.
Remember that $\mathcal{X}_{1,1}$ is characterized by the open conditions
\begin{align*}
x_{12}x_{23}x_{34}x_{14}\ne 0, \qquad  t_{13}\ne t_{24},  \qquad  x_{12}\ne x_{23},
\end{align*}
so $\mathcal{C}\subseteq\mathcal{X}_{1,1}$.
By a result of Thurston \cite{Th80}, $\dim\mathcal{C}\ge 1$. Hence $\deg{\rm gcd}(g,h)>0$.

On the other hand, when $\alpha=1$, so that $\mu^2-\mu+2=0$, we can easily find $\nu=2\mu-3$, $t_{13}=2-2\mu$, $t_{24}=2\mu-4$, $s_{123}=t(\mu-1)^2\ne 0$, and then
\begin{align*}
&t_{24}-2-\frac{(s_{12}+s_{23})^2}{t_{13}-2}-\frac{(s_{12}-s_{23})^2}{\alpha-1-t_{13}}=2\mu-6+2(s_{12}^2+s_{23}^2)   \\
=\ &2(\mu-3+x_{12}^2+x_{23}^2)=2(\mu-3+\mu^2-2\nu)=2-4\mu\ne 0.
\end{align*}
Consequently, $h(\mu)\ne 0$ when $\alpha=1$. By continuity, $h(\mu)\ne 0$ for infinitely many values of $\alpha$.
This contradicts $\deg{\rm gcd}(g,h)>0$.
\end{proof}

Consequently, $(\mu-2\alpha)(\mu+\alpha^2-2\alpha)=-\alpha^3$,
which implies
\begin{align}
\nu=\frac{\alpha\mu}{2\alpha-\mu}=\frac{2}{\alpha}\mu+\alpha-4.   \label{eq:nu}
\end{align}
Therefore, $t_{13}=1$, $t_{24}=\tau$, with
\begin{align*}
\tau=\big(\alpha^2-3\alpha+\frac{4}{\alpha}\big)\mu-\alpha^3+4\alpha^2-7;
\end{align*}
by (\ref{eq:1-t123}), $t_{123}=t_{134}=t\mu/\alpha$. From (\ref{eq:t124+t234}), (\ref{eq:t124-minus-t234}) we solve
$t_{124}=t\varsigma_1$, $t_{234}=t\varsigma_2$, with
\begin{align*}
\varsigma_1&=\big(\frac{\alpha}{2}+\frac{1}{\alpha}-1\big)\mu-\frac{\alpha^2}{2}+2\alpha-2+\frac{\mu}{2\alpha}(t_{12}-t_{23}),  \\
\varsigma_2&=\big(\frac{\alpha}{2}+\frac{1}{\alpha}-1\big)\mu-\frac{\alpha^2}{2}+2\alpha-2-\frac{\mu}{2\alpha}(t_{12}-t_{23}).
\end{align*}

\begin{ass}\label{ass:constraint}
$t_{13}\ne t_{24}$, $x_{12}\ne x_{23}$ if and only if
\begin{align}
\alpha\ne 0,2,4, \qquad (\alpha,\mu)\ne(2+2\sqrt{2},2\sqrt{2}),\ (2-2\sqrt{2},-2\sqrt{2}).  \label{eq:constraint}
\end{align}
\end{ass}

\begin{proof}
As is easy to see, $t_{13}=t_{24}$ is equivalent to
$$\mu=\frac{\alpha^3-4\alpha^2+8}{\alpha^2-3\alpha+4/\alpha};$$
when this holds, (\ref{eq:component-1}) becomes
\begin{align*}
0&=(\alpha^3-4\alpha^2+8)^2+(\alpha^2-4\alpha)(\alpha^3-4\alpha^2+8)\big(\alpha^2-3\alpha+\frac{4}{\alpha}\big)  \\
&\ \ \ \ \ \ +(4\alpha^2-\alpha^3)\big(\alpha^2-3\alpha+\frac{4}{\alpha}\big)^2
=\alpha(\alpha-2)^2(\alpha^2-4\alpha-4).
\end{align*}
Note that if $\alpha^2-4\alpha-4=0$, i.e. $\alpha=2+2\epsilon\sqrt{2}$ with $\epsilon\in\{\pm1\}$, then $\mu=2\epsilon\sqrt{2}$.

Hence $t_{13}\ne t_{24}$ is equivalent to $\alpha\notin\{0,2\}$ and $(\alpha,\mu)\ne(2\pm 2\sqrt{2},\pm2\sqrt{2})$.

The condition $x_{12}=x_{23}$ is equivalent to $\mu^2=4\nu$, i.e. $$0=(4\alpha-\alpha^2)\mu+\alpha^3-4\alpha^2-4\big(\frac{2}{\alpha}\mu+\alpha-4\big)
=\big(4\alpha-\alpha^2-\frac{8}{\alpha}\big)\mu+\alpha^3-4\alpha^2-4\alpha+16;$$
when this holds, (\ref{eq:component-1}) becomes
\begin{align*}
0&=(\alpha^3-4\alpha^2-4\alpha+16)^2+(\alpha^2-4\alpha)(\alpha^3-4\alpha^2-4\alpha+16)\big(\alpha^2-4\alpha+\frac{8}{\alpha}\big)  \\
&\ \ \ \ \ \ +(4\alpha^2-\alpha^3)\big(\alpha^2-4\alpha+\frac{8}{\alpha}\big)^2
=\alpha(\alpha-2)^2(\alpha-4)(\alpha^2-4\alpha-4).
\end{align*}
When $\alpha^2-4\alpha-4=0$, i.e. $\alpha=2+2\epsilon\sqrt{2}$ with $\epsilon\in\{\pm1\}$, we have $\mu=2\epsilon\sqrt{2}$.

So $t_{12}\ne t_{34}$ is equivalent to $\alpha\notin\{0,2,4\}$ and $(\alpha,\mu)\ne(2\pm 2\sqrt{2},\pm2\sqrt{2})$.
\end{proof}

\begin{rmk}
\rm Remember $\alpha\ne 0$, so $\mu\ne\alpha$. Hence $s_{123}\ne 0$ as long as $t\ne 0$.

It is easy to see that $\nu\ne0$, so $x_{12}x_{23}\ne 0$ is ensured.

Also ensured is $t_{12},t_{23}\ne 2$, i.e. $x_{12},x_{23}\ne 1$. Indeed, if $1$ is a root of $x^2-\mu x+\nu$, then
$\mu=(\alpha^2-3\alpha)/(\alpha-2)$, which would contradict (\ref{eq:component-1}).
\end{rmk}

\begin{ass}
The definitive relations all hold.
\end{ass}

\begin{proof}
We have
\begin{align*}
s_{12}+s_{23}&=\mu+1-\alpha,  \\
s_{12}s_{23}&=\nu+\frac{1-\alpha}{2}\mu+\frac{1}{4}(1-\alpha)^2
=\big(\frac{2}{\alpha}+\frac{1-\alpha}{2}\big)\mu+\frac{1}{4}(1+\alpha)^2-4,   \\
s_{123}=s_{134}&=t\big(\frac{1}{\alpha}-\frac{1}{2}\big)(\mu-\alpha),
\end{align*}
and by (\ref{eq:t124+t234}),
$$s_{124}+s_{234}=t\Big((4\alpha-3-\alpha^2-\frac{2}{\alpha})\mu+\alpha^3-5\alpha^2+5\alpha+2\Big).$$

Then (\ref{eq:typeII-1+3}) follows from
\begin{align*}
(\mu+1-\alpha)(\alpha-2)(1-\frac{\mu}{\alpha})+(4\alpha-3-\alpha^2-\frac{2}{\alpha})\mu+\alpha^3-5\alpha^2+5\alpha+2=0,
\end{align*}
which is easy to verify.

If $t=0$, then $s_{ijk}$'s all vanish, so (\ref{eq:typeII-2}) holds. If $t\ne 0$, then (\ref{eq:typeII-2}) is equivalent to (\ref{eq:II-2'}), and
we can verify (\ref{eq:II-2'}) by writing its left-hand-side as
\begin{align*}
&t_{24}-2-\frac{4}{t_{13}+1-\alpha}s_{12}s_{23}-\frac{(\alpha-3)(s_{12}+s_{23})^2}{(t_{13}-2)(\alpha-1-t_{13})}  \\
=\ &\big(\alpha^2-3\alpha+\frac{4}{\alpha}\big)\mu-\alpha^3+4\alpha^2-9
-\frac{4}{2-\alpha}\Big(\big(\frac{2}{\alpha}+\frac{1-\alpha}{2}\big)\mu+\frac{1}{4}(1+\alpha)^2-4\Big)   \\
&-\frac{\alpha-3}{2-\alpha}\big((2+2\alpha-\alpha^2)\mu+\alpha^3-3\alpha^2-2\alpha+1\big)
\end{align*}
and checking that the coefficients of $\mu$ and the constants both cancel out.

We can verify $2s_{123}^2+\det S_{123}=0$ by writing
\begin{align*}
&2s_{123}^2+s_0(s_0+s_{13})(s_0-s_{13})-s_0(s_{12}+s_{23})^2+2(s_0+s_{13})s_{12}s_{23}   \\
=\ &2(\alpha+1)\big(\frac{1}{\alpha}-\frac{1}{2}\big)^2(\mu-\alpha)^2+\frac{(\alpha-3)(2-\alpha)}{2}+\frac{3-\alpha}{2}(\mu+1-\alpha)^2    \\
&-2\Big(\big(\frac{2}{\alpha}+\frac{1-\alpha}{2}\big)\mu+\frac{1}{4}(1+\alpha)^2-4\Big)  \\
=\ &\frac{(\alpha+1)(2-\alpha)^2}{2\alpha^2}\big((2\alpha-\alpha^2)\mu+\alpha^3-3\alpha^2\big)
+\frac{\alpha-3}{2}(2-\alpha)-(1+\alpha)^2+16 \\
&+\frac{3-\alpha}{2}\big((2+2\alpha-\alpha^2)\mu+\alpha^3-3\alpha^2-2\alpha+1\big)
-\big(\frac{4}{\alpha}+1-\alpha\big)\mu
\end{align*}
and checking that the coefficients of $\mu$ and the constants both cancel out.
\end{proof}

\begin{rmk}\label{rmk:no-worry}
\rm When $t=0$, we have $t_{24}=-2$, $s_{123}=s_{124}=s_{134}=s_{234}=0$, $\mu^2+5\mu+5=0$, and $x_1,x_2$ are roots of $x^2-\mu x-2\mu-5=0$.
By Lemma \ref{lem:a1=a3}, $\mathbf{x}_4=\mathbf{x}_2$. It is easy to see $x_1,x_2\ne 1,-3$, so $t_{12},t_{23}\ne 2,t^2-2$.

Recall (\ref{eq:nu}), so that $x_{12},x_{23}$ are the two roots of
\begin{align}
x^2-\mu x+\frac{2}{\alpha}\mu+\alpha-4=0.   \label{eq:root-1}
\end{align}
Each trace is continuous on the algebraic curve determined by (\ref{eq:component-1}),(\ref{eq:constraint}),(\ref{eq:root-1}).

Remember that we have shown ${\rm tr}((\mathbf{r}_1-\mathbf{r}_2)\mathbf{a})=0$ for
$\mathbf{a}\in\{\mathbf{e},\mathbf{x}_1,\mathbf{x}_2,\mathbf{x}_1\mathbf{x}_2\}$. If $t_{12}\ne t^2-2$, then by Remark \ref{rmk:trace-equation} (iii), $\mathbf{r}_1=\mathbf{r}_2$, which implies ${\rm tr}((\mathbf{r}_1-\mathbf{r}_2)\mathbf{x}_3)=0$. When $t_{12}=t^2-2$, by continuity,
${\rm tr}((\mathbf{r}_1-\mathbf{r}_2)\mathbf{x}_3)=0$ still holds; in this case, $s_{123}\ne 0$, implying $\mathbb{C}\langle\mathbf{e},\mathbf{x}_1,\mathbf{x}_2,\mathbf{x}_3\rangle=\mathcal{M}$, so $\mathbf{r}_1=\mathbf{r}_2$ is still ensured.

Similar arguments can be applied to deduce $\mathbf{r}_2=\mathbf{r}_3$ and $\mathbf{r}_3=\mathbf{r}_4$.
\end{rmk}

Writing $t_{12}=x+1$, and replacing $\varsigma_1,\varsigma_2$ respectively by
\begin{align*}
\sigma_+=\big(\frac{1}{\alpha}+1\big)\mu-2+\frac{\mu}{\alpha}x, \qquad  \sigma_-=\frac{\mu^2}{\alpha}+\big(\frac{1}{\alpha}+1\big)\mu-2-\frac{\mu}{\alpha}x,
\end{align*}
we can parameterize $\mathcal{X}_{1,1}$ as
\begin{align*}
\mathcal{X}_{1,1}=\big\{\vec{\mathsf{t}}\colon &t_{12}=t_{14}=x+1, \ t_{23}=t_{34}=\mu+1-x,\ t_{13}=1,\ t_{24}=\tau, \\
&t_{123}=t_{134}=t\mu/\alpha,\ t_{124}=t\sigma_+,\ t_{234}=t\sigma_-,
\ (\ref{eq:component-1}),(\ref{eq:constraint}),(\ref{eq:root-1})\ \text{hold}\big\}.
\end{align*}

Recall that $\mathcal{X}_{1,2}=f^\ast(\mathcal{X}_{1,1})$, so
\begin{align*}
\mathcal{X}_{1,2}=\big\{\vec{\mathsf{t}}\colon &t_{12}=t_{23}=x+1, \ t_{14}=t_{34}=\mu+1-x,\ t_{24}=1,\ t_{13}=\tau, \\
&t_{124}=t_{234}=t\mu/\alpha,\ t_{123}=t\sigma_+,\ t_{134}=t\sigma_-,
\ (\ref{eq:component-1}),(\ref{eq:constraint}),(\ref{eq:root-1})\ \text{hold}\big\}.
\end{align*}

\subsubsection{$t_{13}=t_{24}$, and exactly one of $x_{14}=x_{23}$, $x_{12}=x_{34}$ holds}

Assume $x_{14}=x_{23}$ and $x_{12}\ne x_{34}$. At the end of this subsection we will state the result for the other case, by taking a quarter turn (i.e. acting via $f^\ast$).

Let $z=x_{14}=x_{23}$. Then (\ref{eq:1-2-1})--(\ref{eq:1-2-4}) are equivalent to
$$z^2-\alpha z+\alpha=0.$$
By (\ref{eq:1-1-1}), $t_{123}=t_{124}=:\eta_1$; by (\ref{eq:1-1-3}), $t_{134}=t_{234}=:\eta_2$.

Let
$\mu=x_{12}+x_{34}$, $\nu=x_{12}x_{34}$, and
$$\gamma_1=x_{12}+2-\alpha, \qquad \gamma_2=x_{34}+2-\alpha.$$

The equations (\ref{eq:1-1-2}), (\ref{eq:1-1-4}) are equivalent to
\begin{align}
\eta_1-\eta_2=t\big(1-\frac{1}{z}\big)(x_{12}-x_{34})=\frac{t}{\alpha}z(x_{12}-x_{34}),  \label{eq:remaining2-1}
\end{align}
and (\ref{eq:1-3-1})--(\ref{eq:1-3-4}) are equivalent to
\begin{align}
t(z+2-\alpha)(\eta_1+\eta_2)=2(z+1-\alpha)t_{13}+(2z-\alpha+1)(\mu+2),  \label{eq:remaining2-2}  \\
\gamma_1t\eta_1=(\gamma_1-1)t_{13}+(\gamma_1^2+(\alpha-2)\gamma_1+1)(z+1),  \label{eq:remaining2-3}  \\
\gamma_2t\eta_2=(\gamma_2-1)t_{13}+(\gamma_2^2+(\alpha-2)\gamma_2+1)(z+1).  \label{eq:remaining2-4}
\end{align}

The type II relations (\ref{eq:II-1})--(\ref{eq:II-4}) are equivalent to
\begin{align}
(s_{23}-s_{13})s_{123}&=(s_0-s_{12})s_{234},  \label{eq:II-2-1}  \\
(s_{13}-s_{23})s_{234}&=(s_{34}-s_0)s_{123}.  \label{eq:II-2-2}
\end{align}

\begin{ass}
Either $t=s_{123}=s_{234}=0$, or $ts_{123}s_{234}\ne 0$.
\end{ass}

\begin{proof}
If $t=0$, then by (\ref{eq:remaining2-1}), $\eta_1=\eta_2$, so $s_{234}=s_{123}$. The sum of (\ref{eq:II-2-1}) and (\ref{eq:II-2-2}) implies
$(s_{34}-s_{12})s_{234}=0$; since $x_{34}\ne x_{12}$, we have $s_{123}=s_{234}=0$.

Now suppose $t\ne 0$. We show $s_{123}s_{234}\ne 0$.

If $s_{123}=0\ne s_{234}$, then by (\ref{eq:II-2-1}), $s_{12}=s_0$, so that $t_{12}=t^2-2$. By Lemma \ref{lem:a1=a3},
$\mathbf{x}_1\mathbf{x}_2=\mathbf{x}_2\mathbf{x}_1$, and then by Lemma \ref{lem:priori}, $t_{23}=1$,
contradicting the assumption $x_{23}\ne 0$. Hence $s_{123}=0\ne s_{234}$ is impossible.

Similarly, $s_{123}\ne 0=s_{234}$ is neither possible.

If $s_{123}=s_{234}=0$, then
$$0=s_{123}-s_{234}=\eta_1-\eta_2-\frac{t}{2}(t_{12}-t_{34})=t\big(\frac{z}{\alpha}-\frac{1}{2}\big)(t_{12}-t_{34}),$$
implying $z=\alpha/2$. This together with $z^2=\alpha(z-1)$ implies $\alpha=4$. Now $t_{23}=3=t^2-2$. By Lemma \ref{lem:a1=a3}, $\mathbf{x}_2\mathbf{x}_3=\mathbf{x}_3\mathbf{x}_2$, and then by Lemma \ref{lem:priori}, $t_{34}=1$, contradicting the assumption $x_{34}\ne 0$.
\end{proof}

Suppose $t=0$, i.e., $\alpha=-1$. Then $z^2+z-1=0$. By the above assertion, $\eta_1=\eta_2=0$.
The equation (\ref{eq:remaining2-2}) becomes
\begin{align}
(z+2)t_{13}+(z+1)(\mu+2)=0.   \label{eq:z-mu}
\end{align}
Taking difference between (\ref{eq:remaining2-3}) and (\ref{eq:remaining2-4}) and dividing by $\gamma_1-\gamma_2$, we obtain $t_{13}+(\mu+3)(z+1)=0.$
From this and (\ref{eq:z-mu}) we solve
$$t_{13}=1, \qquad  \mu=-\frac{3z+4}{z+1}=-(3z+4)z=-z-3.$$
Consequently, $\gamma_1+\gamma_2=\mu+6=3-z$.
Summing (\ref{eq:remaining2-3}) and (\ref{eq:remaining2-4}) leads to
\begin{align*}
0=\mu+4+((\mu+6)^2-3\nu-16-2\gamma_1\gamma_2)(z+1)=2z-2(z+1)\nu,
\end{align*}
implying $\nu=z^2$.

\bigskip

From now on, suppose $t\ne 0$. The product of (\ref{eq:II-2-1}) and (\ref{eq:II-2-2}) implies
\begin{align}
(x_{13}-x_{23})^2=(s_{12}-s_0)(s_{34}-s_0)=\gamma_1\gamma_2.   \label{eq:II-2-3}
\end{align}

\begin{ass}\label{ass:ne}
$x_{12},x_{23},x_{34}\ne\alpha-2$.
\end{ass}

\begin{proof}
If $x_{12}=\alpha-2$, then $s_{12}=s_0$; by (\ref{eq:II-2-1}), $s_{13}=s_{23}$, so by (\ref{eq:II-2-2}), $s_{34}=s_0$, contradicting
$t_{34}\ne t_{12}$. Hence $x_{12}\ne\alpha-2$.
Similarly, $x_{34}\ne \alpha-2$.

If $x_{23}=\alpha-2$, then by (\ref{eq:remaining2-2}), $x_{13}=\mu/2$, which together with (\ref{eq:II-2-3})
would imply $\mu^2=4\nu$, contradicting the assumption $x_{12}\ne x_{34}$. Thus, $x_{23}\ne\alpha-2$.
\end{proof}

Rewrite (\ref{eq:remaining2-3}), (\ref{eq:remaining2-4}) as
\begin{align}
t\eta_1&=\Big(1-\frac{1}{\gamma_1}\Big)t_{13}+\Big(x_{12}+\frac{1}{\gamma_1}\Big)(z+1),  \label{eq:remaining2-3'} \\
t\eta_2&=\Big(1-\frac{1}{\gamma_2}\Big)t_{13}+\Big(x_{34}+\frac{1}{\gamma_2}\Big)(z+1).  \label{eq:remaining2-4'}
\end{align}
Taking their difference and using (\ref{eq:remaining2-1}), we obtain
\begin{align*}
\frac{\alpha+1}{\alpha}z(x_{12}-x_{34})=t(\eta_1-\eta_2)=\frac{t_{13}+(\gamma_1\gamma_2-1)(z+1)}{\gamma_1\gamma_2}(x_{12}-x_{34}).
\end{align*}
Hence
$$x_{13}-z=\big(\frac{z}{\alpha}-1\big)\gamma_1\gamma_2=-\frac{1}{z}(x_{13}-z)^2.$$
By (\ref{eq:II-2-3}) and the above assertion, $x_{13}\ne z$. Thus,
$$x_{13}=0.$$

Now $\gamma_1\gamma_2=z^2=\alpha(z-1)$, and the sum of (\ref{eq:remaining2-3'}) and (\ref{eq:remaining2-4'}) gives rise to
\begin{align}
t(\eta_1+\eta_2)&=2-\frac{\mu+4-2\alpha}{\gamma_1\gamma_2}+\Big(\mu+\frac{\mu+4-2\alpha}{\gamma_1\gamma_2}\Big)(z+1)  \nonumber   \\
&=2+\big(z+1+\frac{1}{z}\big)\mu+\frac{4-2\alpha}{z}.      \label{eq:sum-of-eta}
\end{align}
Substituting it into (\ref{eq:remaining2-2}), we obtain
\begin{align*}
0&=(z+2-\alpha)\Big(2+\big(z+1+\frac{1}{z}\big)\mu+\frac{4-2\alpha}{z}\Big) \\
&\ \ \ \ -2(z+1-\alpha)-(2z-\alpha+1)(\mu+2)  \\
&=\Big(z+2-\alpha+\frac{2-\alpha}{z}\Big)\mu+4-4z+\frac{2(2-\alpha)^2}{z}.
\end{align*}
Hence
\begin{align}
\mu&=\frac{z(4z-4)-2(2-\alpha)^2}{z(z+2-\alpha)+2-\alpha}
=\frac{(2\alpha-2)z-\alpha^2+2\alpha-4}{z+1-\alpha} \nonumber  \\
&=(1-z)\big((2\alpha-2)z-\alpha^2+2\alpha-4\big)=2z+(2-\alpha)z^2+2\alpha-4.  \label{eq:mu}
\end{align}
Using this (and $z^2=\alpha z-\alpha$), we can rewrite (\ref{eq:sum-of-eta}) as
\begin{align*}
t(\eta_1+\eta_2)=2+(z+1)\mu+2+(2-\alpha)z=\alpha(\alpha+1)((3-\alpha)z+\alpha-2).
\end{align*}
Thus, $\eta_1+\eta_2=t\alpha((3-\alpha)z+\alpha-2)$. From this and (\ref{eq:remaining2-1}) we can solve
\begin{align}
\eta_1=t\big(1+\frac{z}{\alpha}\gamma_1\big),   \qquad   \eta_2=t\big(1+\frac{z}{\alpha}\gamma_2\big).   \label{eq:eta}
\end{align}

\begin{ass}
The definitive relations all hold.
\end{ass}

\begin{proof}
It follows from (\ref{eq:eta}) that
$$s_{123}=t\Big((\frac{z}{\alpha}-\frac{1}{2})\gamma_1+1-\frac{z}{2}\Big),  \qquad
s_{234}=t\Big((\frac{z}{\alpha}-\frac{1}{2})\gamma_2+1-\frac{z}{2}\Big).$$

Note that $\gamma_i^2=(\mu+4-2\alpha)\gamma_i-\alpha(z-1)$, and
\begin{align*}
\det S_{123}&=s_0^3-2(s_0+1)(\gamma_1+s_0)(z-1-s_0) \\
&\ \ \ \ -s_0((\gamma_1+s_0)^2+(s_0+1)^2+(z-1-s_0)^2)  \\
&=2(2s_0+1)\gamma_1-2(s_0+1)z\gamma_1-s_0(\gamma_1^2+z^2).
\end{align*}
Then it is easy to verify $2s_{123}^2+\det S_{123}=0$.

The type II relations (\ref{eq:II-2-1}), (\ref{eq:II-2-2}) amount to $zs_{123}=-\gamma_1s_{234}$, which can be verified by checking
\begin{align*}
&z\Big((\frac{z}{\alpha}-\frac{1}{2})\gamma_1+1-\frac{z}{2}\Big)+\gamma_1\Big((\frac{z}{\alpha}-\frac{1}{2})\gamma_2+1-\frac{z}{2}\Big)  \\
=\ &\big(\frac{z^2}{\alpha}+1-z\big)\gamma_1+z-\frac{z^2}{2}+\big(\frac{z}{\alpha}-\frac{1}{2}\big)\alpha(z-1)=0.
\end{align*}

Therefore, the definitive relations all hold.
\end{proof}

\begin{ass}
$x_{12}\ne x_{34}$ if and only if
$$\alpha\ne 0,2,4, \qquad (\alpha,z)\ne(2+2\sqrt{2},\sqrt{2}),\ (2-2\sqrt{2},-\sqrt{2}).$$
\end{ass}

\begin{proof}
Clearly, $x_{12}=x_{34}$, which is the same as $\gamma_1=\gamma_2$, is equivalent to
\begin{align*}
0=(\mu+4-2\alpha)^2-4z^2=(\alpha^2-2\alpha)(z-1)\big((\alpha^2-2\alpha-4)z-\alpha^2+2\alpha\big).
\end{align*}
With $z^2-\alpha z+\alpha=0$ recalled, $(\alpha^2-2\alpha-4)z-\alpha^2+2\alpha=0$ is equivalent to
$$\alpha(\alpha-4)(\alpha^2-4\alpha-4)=0.$$

If $\alpha^2-4\alpha-4=0$, i.e. $\alpha=2+2\epsilon\sqrt{2}$ with $\epsilon\in\{\pm1\}$, then we can compute
$$z=\frac{\alpha^2-2\alpha}{\alpha^2-2\alpha-4}=\epsilon\sqrt{2},$$

Thus, $x_{12}\ne x_{34}$ if and only if
$\alpha\ne 0,2,4$ and $(\alpha,z)\ne(2\pm2\sqrt{2},\pm\sqrt{2})$.
\end{proof}

In view of (\ref{eq:mu}), $\gamma_1,\gamma_2$ are the two roots of
$$\gamma^2+((\alpha-2)z^2-2z)\gamma+z^2=0.$$

\begin{rmk}\label{rmk:allow}
\rm When $t_{12}=1$, the whole procedure still works. In this case, (\ref{eq:remaining2-3}) is redundant, but it makes the equations symmetric.

On the other hand, the case $t_{12}=t_{34}=1$ has been automatically excluded: if $\gamma_1=\gamma_2=2-\alpha$, then $(2-\alpha)^2=z^2$
and $4-2\alpha=2z-(\alpha-2)z^2$, which contradicts $z^2-\alpha z+\alpha=0$.

Furthermore, by Assertion \ref{ass:ne}, $t_{12},t_{23},t_{34}\ne t^2-2$. Clearly $z\ne 1$, so that $t_{23}=t_{14}\ne 2$. When $t_{12}=2$ (resp. $t_{34}=2$), similarly as Remark \ref{rmk:no-worry}, we can still ensure $\mathbf{r}_1=\mathbf{r}_2$ (resp. $\mathbf{r}_3=\mathbf{r}_4$).
\end{rmk}

Unifying the cases $t=0$ and $t\ne0$, we obtain a component
\begin{align*}
\mathcal{X}_{2,1}=\big\{\vec{\mathsf{t}}\colon &t_{12}=\gamma+\alpha-1,\ t_{34}=z^2/\gamma+\alpha-1,\ t_{23}=t_{14}=z+1, \\
&t_{13}=t_{24}=1,\ t_{123}=t_{124}=t+tz\gamma/\alpha,\ t_{134}=t_{234}=t+tz^3/(\alpha\gamma), \\
&z^2-\alpha z+\alpha=0,\ \ \gamma^2+((\alpha-2)z^2-2z)\gamma+z^2=0, \\
&\alpha\ne 0,2,4, \ \  (\alpha,z)\ne(2+2\sqrt{2},\sqrt{2}),(2-2\sqrt{2},-\sqrt{2})\big\}.
\end{align*}

Taking a quarter turn leads to
\begin{align*}
\mathcal{X}_{2,2}=\big\{\vec{\mathsf{t}}\colon &t_{23}=\gamma+\alpha-1,\ t_{14}=z^2/\gamma+\alpha-1,\ t_{12}=t_{34}=z+1, \\
&t_{13}=t_{24}=1,\ t_{134}=t_{234}=t+tz\gamma/\alpha,\ t_{123}=t_{124}=t+tz^3/(\alpha\gamma),\\
&z^2-\alpha z+\alpha=0,\ \ \gamma^2+((\alpha-2)z^2-2z)\gamma+z^2=0, \\
&\alpha\ne 0,2,4, \ \  (\alpha,z)\ne(2+2\sqrt{2},\sqrt{2}),(2-2\sqrt{2},-\sqrt{2})\big\}.
\end{align*}

\subsubsection{$t_{13}=t_{24}$, $x_{14}=x_{23}$, $x_{12}=x_{34}$}

As is easy to see, (\ref{eq:1-1-1})--(\ref{eq:1-1-4}) are equivalent to
$t_{123}=t_{124}=t_{134}=t_{234}=:\eta$, so $s_{123}=s_{124}=s_{134}=s_{234}=:\tilde{\eta}.$

Obviously, (\ref{eq:1-2-1})--(\ref{eq:1-2-4}) hold, and (\ref{eq:1-3-1})--(\ref{eq:1-3-4}) are equivalent to
\begin{align}
t(\alpha-1-t_{12})\eta+(t_{12}-\alpha)t_{13}=(\alpha t_{12}-t_{12}^2-\alpha)t_{23},  \label{eq:trace3-1}   \\
t(\alpha-1-t_{23})\eta+(t_{23}-\alpha)t_{13}=(\alpha t_{23}-t_{23}^2-\alpha)t_{12},  \label{eq:trace3-2}
\end{align}
Each of (\ref{eq:II-1})--(\ref{eq:II-4}) is equivalent to
\begin{align}
(s_{12}+s_{23}-s_0-s_{13})\tilde{\eta}=0.   \label{eq:II-equivalent}
\end{align}

Observe that
\begin{align}
\det S_\diamond&=\big((s_{13}+s_0)^2-(s_{12}+s_{23})^2\big)\big((s_{13}-s_0)^2-(s_{12}-s_{23})^2\big), \label{eq:det-3} \\
\det S_{123}&=s_0^3+2(s_{13}\pm s_0)s_{12}s_{23}-s_0s_{13}^2-s_0(s_{12}\pm s_{23})^2.  \label{eq:det-4}
\end{align}

\begin{ass}
$x_{12}=x_{23}$.
\end{ass}

\begin{proof}
Assume on the contrary that $x_{12}\ne x_{23}$.

Let $\mu=x_{12}+x_{23}$, $\nu=x_{12}x_{23}$.
From (\ref{eq:trace3-1}), (\ref{eq:trace3-2}) we can solve
\begin{align}
t\eta&=\alpha\mu-\alpha^2+2\alpha,   \label{eq:teta}   \\
t_{13}&=(\alpha-1)\mu-\alpha^2+3\alpha-1-\nu.  \label{eq:t13}
\end{align}

If $\tilde{\eta}\ne 0$, then $s_{12}+s_{23}=s_0+s_{13}$, implying $t_{13}=\mu+3-\alpha$. It follows from (\ref{eq:t13}) that
$\nu=(\alpha-2)(\mu+2-\alpha)$, and we can compute
\begin{align*}
s_{12}s_{23}-s_0s_{13}=\nu+\frac{1-\alpha}{2}\mu+\frac{1}{4}(1-\alpha)^2+\frac{3-\alpha}{2}\big(\mu+\frac{5-3\alpha}{2}\big)=0.
\end{align*}
From (\ref{eq:det-4}) we see
$$\det S_{123}=(s_{13}+s_0)(s_{12}s_{23}-s_0s_{13})=0,$$
contradicting $2s_{123}^2+\det S_{123}=0$.
Hence $\tilde{\eta}=0$.

Then $2\eta=t(2+\mu+t_{13})-t^3$,
which together with (\ref{eq:teta}), (\ref{eq:t13}) imply
\begin{align}
(\alpha+1)\nu=\alpha(\alpha-1)\mu+3\alpha^2-2\alpha-\alpha^3.   \label{eq:multiple-of-nu}
\end{align}

It follows from $\det S_\diamond=0$ that $(s_{13}+\epsilon s_0)^2=(s_{12}+\epsilon s_{23})^2$ for some $\epsilon\in\{\pm1\}$.
Hence $0=2s_{123}^2+\det S_{123}=\det S_{123}$ becomes
\begin{align*}
0&=s_0^3+2(s_{13}+\epsilon s_0)s_{12}s_{23}-s_0s_{13}^2-s_0(s_{13}+\epsilon s_0)^2  \\
&=2(s_{13}+\epsilon s_0)(s_{12}s_{23}-s_0s_{13}).
\end{align*}

If $s_{12}s_{23}=s_0s_{13}$, which reads
$$\nu+\frac{1-\alpha}{2}\mu+\frac{(1-\alpha)^2}{4}=\frac{\alpha-3}{2}\Big((\alpha-1)\mu-\alpha^2+3\alpha-1-\nu-\frac{\alpha+1}{2}\Big),$$
we would deduce
$$\alpha(\alpha-1)(\nu+(\alpha-2)^2)=\alpha(\alpha-1)(\alpha-2)\mu\stackrel{(\ref{eq:multiple-of-nu})}
=(\alpha-2)\big((\alpha+1)\nu+\alpha^3-3\alpha^2+2\alpha\big),$$
which contradicts $\nu\ne0$. Hence $s_{12}s_{23}\ne s_0s_{13}$.

Then $s_{13}+\epsilon s_0=0$, so that $s_{12}+\epsilon s_{23}=0$, forcing $\epsilon=1$.
Consequently, $0=s_{13}+s_0=t_{13}-2$, implying $t_{13}=2$, and $s_{12}+s_{23}=0$ implies $\mu=\alpha-1$, then by (\ref{eq:t13}), $\nu=\alpha-2$. However, this violates (\ref{eq:multiple-of-nu}). Thus, $x_{12}=x_{23}$.
\end{proof}

Now (\ref{eq:trace3-1}) becomes
\begin{align}
t(\alpha-1-t_{12})\eta+(t_{12}-\alpha)t_{13}+(t_{12}^2-\alpha t_{12}+\alpha)t_{12}=0,  \label{eq:trace'}
\end{align}
and (\ref{eq:II-equivalent}) becomes
$2s_{12}\tilde{\eta}=(s_{13}+s_0)\tilde{\eta}.$
Moreover, (\ref{eq:det-3}), (\ref{eq:det-4}) become
\begin{align*}
\det S_\diamond&=\big((s_{13}+s_0)^2-4s_{12}^2\big)(s_{13}-s_0)^2,  \\
\det S_{123}&=(s_0-s_{13})\big(s_0(s_0+s_{13})-2s_{12}^2\big).
\end{align*}

Suppose $\tilde{\eta}=0$. We have $t_{12}\ne 2,t^2-2$: otherwise by Lemma \ref{lem:a1=a3}, $\mathbf{x}_1\mathbf{x}_2=\mathbf{x}_2\mathbf{x}_1$, and then by Lemma \ref{lem:priori}, $t_{23}=1$, contradicting the assumption $x_{23}\ne 0$. Thus, the irreducibility has been guaranteed.

It follows from $\tilde{\eta}=0$ that $2\eta=t(2t_{12}+t_{13}-t^2)$. From $\det S_\diamond=\det S_{123}=0$ we deduce that
either (i) $s_{13}=s_0$, or (ii) $s_{13}=-s_0$, $s_{12}=0$.

\begin{enumerate}
  \item[\rm(i)] $t_{13}=t^2-2$, $\eta=t(t_{12}-1)$, and (\ref{eq:trace'}) becomes
        $$(t_{12}-t^2+2)\big(t_{12}^2-(t^2+1)t_{12}+2t^2-1\big)=0,$$
        implying
        $$t_{12}^2-(t^2+1)t_{12}+2t^2-1=0.$$
        By Lemma \ref{lem:irreducible}, $\mathbf{C}\langle\mathbf{e},\mathbf{x}_1,\mathbf{x}_2,\mathbf{x}_1\mathbf{x}_2\rangle=\mathcal{M}$.
        Now that ${\rm tr}(\mathbf{x}_1\mathbf{a})={\rm tr}(\mathbf{x}_3\mathbf{a})$ for $\mathbf{a}\in\{\mathbf{e},\mathbf{x}_1,\mathbf{x}_2,\mathbf{x}_1\mathbf{x}_2\}$,
        we have $\mathbf{x}_1=\mathbf{x}_3$.
        Similarly, $\mathbf{x}_4=\mathbf{x}_2$. As is easy to check, $t_{12}\notin\{1,2,t^2-2\}$ is equivalent to $t^2\ne 1,5$.
  \item[\rm(ii)] $t_{13}=2$, $t_{12}=t^2/2$, $\eta=t$, so (\ref{eq:trace'}) becomes
        $(t^2-4)(t^4-6t^2+4)=0$, which together with $t_{12}\ne 2$ implies $t^4-6t^2+4=0$.

        In this case, by Lemma \ref{lem:a1=a3}, $\mathbf{x}_3=\mathbf{x}_1^{-1}$, and also $\mathbf{x}_4=\mathbf{x}_2^{-1}$.
\end{enumerate}

Now suppose $\tilde{\eta}\ne 0$. By (\ref{eq:II-equivalent}), $s_{13}+s_0=2s_{12}$, so that
$t_{13}=2t_{12}+1-\alpha$, and (\ref{eq:trace'}) becomes
$$(t_{12}+1-\alpha)t\eta=(t_{12}+1-\alpha)(t_{12}^2+t_{12}-\alpha).$$
The assumption $s_{123}=\tilde{\eta}\ne 0$ forces
$$0\ne\det S_{123}=4(s_{12}-s_0)^2s_{12}=4(t_{12}+1-\alpha)^2\Big(t_{12}-\frac{t^2}{2}\Big),$$
implying $t_{12}+1-\alpha\ne 0$, i.e. $t_{12}\ne t^2-2$.
Hence $t\eta=t_{12}^2+t_{12}-\alpha$.
Now
$$s_{123}=\eta-\frac{t}{2}(2t_{12}+t_{13}-t^2)=\eta+t(t^2-1-2t_{12}).$$
Then $2s_{123}^2+\det S_{123}=0$ is equivalent to
\begin{align*}
0&=\big(\eta+t(t^2-1-2t_{12})\big)^2+(t_{12}+2-t^2)^2(2t_{12}-t^2)   \\
&=\eta^2+t(t^2+4-2t_{12})\eta-5t^2t_{12}+t^4+4t^2-6;
\end{align*}
in deducing the second line, we have repeatedly used $t_{12}^2=t\eta-t_{12}+\alpha$.

\begin{rmk}
\rm We can include the above Case (ii) by allowing $s_{123}=0$.

In the case $s_{123}\ne 0$, if $t_{12}=2$, then similarly as Remark \ref{rmk:no-worry},
$\mathbf{r}_1=\mathbf{r}_2=\mathbf{r}_3=\mathbf{r}_4$ can still be ensured.
\end{rmk}

Summarizing the result of this subsection, put
\begin{align*}
\mathcal{X}'_3=\big\{\vec{\mathsf{t}}\colon &t_{12}=t_{23}=t_{34}=t_{14}=u,\ t_{123}=t_{124}=t_{134}=t_{234}=tu-t,  \\
&t_{13}=t_{24}=t^2-2, \ u^2-(t^2+1)u+2t^2-1=0,\  t^2\ne 1,5\big\},  \\
\mathcal{X}_4=\big\{\vec{\mathsf{t}}\colon &t_{12}=t_{23}=t_{34}=t_{14}=u,\ t_{13}=t_{24}=2u+2-t^2, \\
&t_{123}=t_{124}=t_{134}=t_{234}=\eta, \ \ t\eta=u^2+u-t^2+1,  \\
&\eta^2+t(t^2+4-2u)\eta-5t^2u+t^4+4t^2-6=0,  \ u\ne 1\big\}.
\end{align*}

\subsection{$x_{12}x_{23}x_{34}x_{14}=0$}

\subsubsection{Exactly one of $x_{12},x_{23},x_{34},x_{14}$ vanishes}

Up to the rotational symmetry, we only need to consider the case $x_{12}=0$ and $x_{23}x_{34}x_{14}\ne 0$. The other three cases are similar.

By Remark \ref{rmk:trace-equation} (ii), $t_{14}=t_{23}$ and
\begin{align*}
(t^2-1)(t_{124}-t_{123})+t(t_{13}-t_{24})=0,  \\
(t^2-2)tt_{124}-tt_{123}+t_{13}+(1-t^2)t_{24}+t_{14}=0.
\end{align*}
Now (\ref{eq:1-1-2})--(\ref{eq:1-1-4}) are valid, implying $t_{123}=t_{124}$ and $t_{134}=t_{234}$;
also, (\ref{eq:1-2-2})--(\ref{eq:1-2-4}) are valid, implying $t_{13}=t_{24}$ and $x_{23}^2-\alpha x_{23}-\alpha=0$.

This part can be merged into $\mathcal{X}_{2,1}$ (given in Section 4.1.2); see Remark \ref{rmk:allow}.

\subsubsection{$x_{12}=x_{34}=0\ne x_{14}x_{23}$ or $x_{14}=x_{23}=0\ne x_{12}x_{34}$}

Similarly as above, we only deal with the case $x_{12}=x_{34}=0\ne x_{14}x_{23}$.

Let $t_{14}=t_{23}=:y\ne 1$.
By Remark \ref{rmk:trace-equation} (i),(iv), we have $t_{13}=t_{24}=:u$, $t_{123}=t_{234}$, and $t_{134}=t_{124}$.
Hence $s_{123}=s_{234}$, and $s_{134}=s_{124}$.

By Remark \ref{rmk:trace-equation} (ii),(iv),
\begin{align}
(t^2-1)(t_{124}-t_{123})=0,  \\
(t^2-2)tt_{124}-tt_{123}+(2-t^2)u+y=0.
\end{align}
Rotating (\ref{eq:tr1+tr2}) by $\pi/4$ and $-\pi/4$ respectively yields
\begin{align*}
t(t^2-2-y)(t_{123}+t_{234})+2(y+1-t^2)u=2(1-y^2+(t^2-1)y-t^2),  \\
t(t^2-2-y)(t_{134}+t_{124})+2(y+1-t^2)u=2(1-y^2+(t^2-1)y-t^2).
\end{align*}
Their difference gives rise to $t(t^2-2-y)(t_{123}-t_{124})=0$.

If $t_{124}\ne t_{123}$, then $t^2=1$ and $y=-1$, implying $u=1$, so $s_{12}=s_{13}=1/2$, $s_0=s_{23}=-3/2$. It is easy to see
$\det S_{123}=0$. Hence $s_{123}=0$. Similarly, $\det S_{124}=0$, implying $s_{124}=0$. But this would imply $t_{124}=t_{123}$.

Thus, $t_{123}=t_{234}=t_{134}=t_{124}=:\eta$, and the trace equations become
\begin{align*}
(t^2-3)t\eta+(2-t^2)u+y&=0,   \\
t(t^2-2-y)\eta+(y+1-t^2)u&=1-y^2+(t^2-1)y-t^2,
\end{align*}
from which we can solve
\begin{align*}
t\eta&=(t^2-1)y+(1-t^2)(t^2-2),   \\
u&=(t^2-2)y+(1-t^2)(t^2-3).
\end{align*}

\begin{ass}
$t^2\ne 3$.
\end{ass}

\begin{proof}
If $t^2=3$, then $t\eta=2(y-1)$, and $u=y$, so
$$ts_{123}=t\eta-\frac{t^2}{2}(u+y+1-t^2)=2(y-1)-\frac{3}{2}(2y-2)=1-y\ne 0.$$
On the other hand, $s_0=-1/2=s_{12}$, and $s_{13}=s_{23}$, so $\det S_{123}=0$. This contradicts $2s_{123}^2+\det S_{123}=0$.
\end{proof}

\medskip

If $t^2=1$, then $\eta=u+y=0$, so $s_{123}=s_{124}=0$, and we have $s_0=-3/2$, $s_{12}=1/2$, $s_{13}+s_{23}=-1$. It follows from
$$0=\det S_{123}=-\frac{27}{8}+2\cdot \frac{1}{2}s_{13}s_{23}+\frac{3}{2}\big(s_{13}^2+s_{23}^2+\frac{1}{4}\big)=-2s_{13}s_{23}-\frac{3}{2}$$
that $\{s_{13},s_{23}\}=\{1/2,-3/2\}$. Since $s_{23}=y-1/2\ne 1/2$, we have $s_{23}=-3/2$, $s_{13}=1/2$.
Hence $t_{13}=t_{24}=1$, $t_{23}=t_{14}=-1=t^2-2$. By Lemma \ref{lem:a1=a3}, $\mathbf{x}_2=\mathbf{x}_3$, $\mathbf{x}_1=\mathbf{x}_4$.
It is easy to check $\mathbf{r}_1=\mathbf{r}_2=\mathbf{r}_3=\mathbf{r}_4$.

\medskip

Now suppose $t^2\ne 1$.
\begin{ass}
$y=t^2-2$.
\end{ass}

\begin{proof}
From (\ref{eq:II-1}) we see $(s_{13}+s_0)s_{123}=(s_{12}+s_{23})s_{123}$, so that
\begin{align*}
0&=(s_{13}+s_0-s_{12}-s_{23})ts_{123}=(u-y+t^2-3)\big(t\eta-\frac{t^2}{2}(u+y+1-t^2)\big)   \\
&=\frac{1}{2}(t^2-3)(2-t^2)(t^2-1)(y+2-t^2)^2.
\end{align*}
Hence $y=t^2-2$ or $t^2=2$.

If $y\ne t^2-2$, then $t^2-2=0\ne y$, so $u=1$, $\eta=y/t$, and
$$s_0=-1, \qquad   s_{12}=s_{13}=0, \qquad  s_{23}=y-1, \qquad  s_{123}=s_{234}=0;$$
it follows from
$0=2s_{123}^2+\det S_{123}=y^2-2y$
that $y=2$. Hence $t^2\ne 4$. By Lemma \ref{lem:a1=a3}, $\mathbf{x}_3=\mathbf{x}_2^{-1}$, then by Lemma \ref{lem:priori}, $\mathbf{x}_3=\mathbf{x}_2$, implying $\mathbf{x}_2^2=\mathbf{e}$; by (\ref{eq:matrix-identity-2-1}), $\mathbf{x}_2\in\{\pm\mathbf{e}\}$.
This contradicts the irreducibility.
\end{proof}

Consequently, $u=1$, $s_{12}=s_{13}$, and $s_{23}=s_0$, implying $\det S_{123}=0$, which forces $s_{123}=0$.
By Lemma \ref{lem:a1=a3}, $\mathbf{x}_2\mathbf{x}_3=\mathbf{x}_3\mathbf{x}_2$. By Lemma \ref{lem:priori}, $\mathbf{x}_2=\mathbf{x}_3$, and $\mathbf{x}_4=\mathbf{x}_1$. From the proof of Lemma \ref{lem:priori} we see  $\mathbf{x}_1\mathbf{x}_2\mathbf{x}_1=\mathbf{x}_2\mathbf{x}_1\mathbf{x}_2$, ensuring $\mathbf{r}_1=\mathbf{r}_2=\mathbf{r}_3=\mathbf{r}_4$.
Furthermore,
$$\eta=t_{123}={\rm tr}(\mathbf{x}_1\mathbf{x}_2^2)={\rm tr}(\mathbf{x}_1(t\mathbf{x}_2-\mathbf{e}))=0.$$

We can unify the cases $t^2=1$ and $t^2\ne 1$, to get a component
\begin{align*}
\mathcal{X}_{5,1}=\big\{\vec{\mathsf{t}}\colon &t_{12}=t_{34}=1,\ t_{23}=t_{14}=t^2-2,\ t_{13}=t_{24}=1, \\
&t_{123}=t_{124}=t_{134}=t_{234}=0, \ t^2\ne 3\big\}.
\end{align*}

Taking a quarter turn, we obtain
\begin{align*}
\mathcal{X}_{5,2}=\big\{\vec{\mathsf{t}}\colon &t_{23}=t_{14}=1,\ t_{12}=t_{34}=t^2-2,\ t_{13}=t_{24}=1, \\
&t_{123}=t_{124}=t_{134}=t_{234}=0, \ t^2\ne 3\big\}.
\end{align*}

\subsubsection{$x_{i-1,i}=x_{i,i+1}=0$ for some $i$}

By Remark \ref{rmk:trace-equation}, $t_{12}=t_{23}=t_{34}=t_{14}=1$.

The remaining trace equations are
\begin{align}
(t^2-1)(t_{124}-t_{123})=t(t_{24}-t_{13}),  \label{eq:remaining4-1}  \\
(t^2-2)tt_{124}-tt_{123}=(t^2-1)t_{24}-t_{13}-1,  \label{eq:remaining4-2}  \\
(t^2-1)(t_{123}-t_{234})=t(t_{13}-t_{24}), \nonumber   \\
(t^2-2)tt_{123}-tt_{234}=(t^2-1)t_{13}-t_{24}-1,  \nonumber  \\
(t^2-1)(t_{234}-t_{134})=t(t_{24}-t_{13}), \nonumber  \\
(t^2-2)tt_{234}-tt_{134}=(t^2-1)t_{24}-t_{13}-1,  \nonumber  \\
(t^2-1)(t_{134}-t_{124})=t(t_{13}-t_{24}), \nonumber \\
(t^2-2)tt_{134}-tt_{124}=(t^2-1)t_{13}-t_{24}-1.  \nonumber
\end{align}
From these we can see (no matter $t^2=1$ or not)
$$t_{234}=t_{124}, \qquad  t_{134}=t_{123},$$
and then the third to eighth equations all follow from (\ref{eq:remaining4-1}), (\ref{eq:remaining4-2}).

The type II relations (\ref{eq:II-1})--(\ref{eq:II-4}) are equivalent to
\begin{align}
(2-t^2)s_{123}=(t_{13}-2)s_{234}, \qquad  (t_{24}-2)s_{123}=(2-t^2)s_{234}.   \label{eq:typeII-4}
\end{align}

By direct computation,
\begin{align}
\det S_{123}&=(s_0-s_{13})\big(s_0^2+s_0s_{13}-\frac{1}{2}(2-t^2)^2\big),  \label{eq:det-1}  \\
\det S_\diamond&=(s_{13}-s_0)(s_{24}-s_0)((s_{13}+s_0)(s_{24}+s_0)-(2-t^2)^2).   \label{eq:det-2}
\end{align}
So $\det S_\diamond=0$ if and only if one of the following holds:
$$t_{13}=t^2-2, \qquad t_{24}=t^2-2,  \qquad  (t_{13}-2)(t_{24}-2)=(2-t^2)^2.$$

\begin{ass}
$t^2\ne3$.
\end{ass}

\begin{proof}
Assume $t^2=3$. Then from (\ref{eq:remaining4-1}), (\ref{eq:remaining4-2}) we see $t_{24}=2-t_{13}$, and from $\det S_\diamond=0$ we deduce $t_{13}=1$.
It follows that $\det S_{123}=0$, implying $s_{123}=0$. Since $t_{12}=t_{13}=t^2-2$,
applying Lemma \ref{lem:a1=a3} we can deduce $\mathbf{x}_1=\mathbf{x}_2=\mathbf{x}_3$, but this contradicts the irreducibility.
\end{proof}

Suppose $t^2=1$. Then (\ref{eq:remaining4-1}), (\ref{eq:remaining4-2}) can be replaced by
$$t_{13}=t_{24},  \qquad  t_{124}+t_{123}=t(t_{13}+1).$$
We have $s_0=-3/2$, $s_{12}=s_{23}=1/2$, and (\ref{eq:typeII-4}) becomes
$$s_{123}=(t_{13}-2)s_{234}, \qquad  s_{234}=(t_{13}-2)s_{123}.$$
\begin{enumerate}
  \item If $s_{123}=0$, then $s_{234}=0$, and from $\det S_{123}=\det S_\diamond=0$ we see
        $t_{13}=-1=t_{24}$. By Lemma \ref{lem:a1=a3}, $\mathbf{x}_3=\mathbf{x}_1$, and $\mathbf{x}_4=\mathbf{x}_2$.
        Then it is easy to verify
        $\mathbf{r}_1=\mathbf{r}_2=\mathbf{r}_3=\mathbf{r}_4$.
  \item If $s_{123}s_{234}\ne 0$, then $(t_{13}-2)^2=1$, which together with $\det S_{123}=-2s_{123}^2\ne0$ implies $t_{13}=1$,
        so $\det S_{123}=-2$. From $2s_{123}^2+\det S_{123}=0$ we deduce $s_{123}=\pm1=\epsilon t$ with $\epsilon\in\{\pm1\}$, so that $t_{123}=s_{123}+t=t(1+\epsilon)$, and $t_{234}=t(1-\epsilon)$.
\end{enumerate}

\medskip

Suppose $t=0$. Then (\ref{eq:remaining4-1}), (\ref{eq:remaining4-2}) are just
$t_{124}=t_{123}$, $t_{13}+t_{24}=-1,$
and (\ref{eq:II-1})--(\ref{eq:II-4}) become
$$2t_{123}=(t_{13}-2)t_{123}, \qquad  2t_{234}=(t_{24}-2)t_{123}.$$
\begin{itemize}
  \item If $t_{13},t_{24}\ne-2$, then $\det S_\diamond=0$ implies
        $4=(t_{13}-2)(t_{24}-2)=t_{13}t_{24}+6,$
        i.e. $t_{13}t_{24}=-2$, which combined with $t_{13}+t_{24}=-1$ yields a contradiction.
  \item If $t_{13}=-2$, then $t_{24}=1$, $t_{123}=0=t_{234}$, and the definitive relations are easy to verify.
  \item Similarly, if $t_{24}=-2$, then $t_{13}=1$, $t_{123}=0=t_{234}$, and the definitive relations are easy to verify.
\end{itemize}

\medskip

From now on, assume $t^2\ne 0,1$. Remember $t^2\ne 3$.

\begin{ass}
$t_{13},t_{24}\ne t^2-2$.
\end{ass}

\begin{proof}
If $t_{13}=t^2-2$, then by (\ref{eq:det-1}), $\det S_{123}=0$, so $s_{123}=0$, which implies $s_{234}=0$, and then
\begin{align*}
t_{123}&=s_{123}+\frac{t}{2}(t_{12}+t_{23}+t_{13}-t^2)=0,  \\
t_{124}=t_{234}&=s_{234}+\frac{t}{2}(t_{23}+t_{34}+t_{24}-t^2)=\frac{t}{2}(t_{24}-t^2+2).
\end{align*}
The equation (\ref{eq:remaining4-1}) becomes
$$(t^2-1)t_{124}=t(t_{24}-t^2+2)=2t_{124},$$
which implies $t_{124}=0$, so that $t_{24}=t^2-2$. But (\ref{eq:remaining4-2}) does not hold, due to $(t^2-2)^2-1\ne 0$.
Thus, $t_{13}\ne t^2-2$.

Similarly, $t_{24}\ne t^2-2$.
\end{proof}

Let $\delta=t^2-2$.
From (\ref{eq:remaining4-1}), (\ref{eq:remaining4-2}) we obtain
\begin{align}
t_{123}&=\frac{(t^4-3t^2+1)t_{13}+t_{24}+1-t^2}{t(t^2-1)(t^2-3)},  \label{eq:t123}  \\
t_{124}&=\frac{(t^4-3t^2+1)t_{24}+t_{13}+1-t^2}{t(t^2-1)(t^2-3)}.  \label{eq:t124}
\end{align}
Then by (\ref{eq:typeII-4}),
\begin{align*}
0&=(t^2-t_{24})s_{123}+(t_{13}-t^2)s_{234}   \\
&=(t^2-t_{24})\big(t_{123}+\frac{t}{2}(\delta-t_{13})\big)+(t_{13}-t^2)\big(t_{234}+\frac{t}{2}(\delta-t_{24})\big)   \\
&=\frac{(t_{13}-t_{24})(t_{13}+t_{24}+\delta^2-3)}{t(t^2-1)(t^2-3)}.
\end{align*}
Hence $t_{13}+t_{24}=3-\delta^2$ or $t_{13}=t_{24}$.

Note that $\det S_\diamond=0$ is equivalent to
$$(t_{13}-2)(t_{24}-2)=\delta^2.$$
So there are four possibilities:
(i) $t_{13}=1$, $t_{24}=2-\delta^2$; (ii) $t_{24}=1$, $t_{13}=2-\delta^2$; (iii) $t_{13}=t_{24}=2+\delta$; (iv) $t_{13}=t_{24}=2+\delta$.
\begin{enumerate}
  \item[\rm(i)] If $t_{13}=1$, $t_{24}=2-\delta^2$, then by (\ref{eq:t123}), $t_{123}=0$, and by (\ref{eq:t124}), $t_{124}=3t-t^3$.
        Hence $s_{123}=t(t^2-3)/2\ne 0$, and we can use (\ref{eq:det-1}) to verify
        $$\det S_{123}=-\frac{1}{2}t^2(t^2-3)^2=-2s_{123}^2.$$
  \item[\rm(ii)] Similarly, if $t_{24}=1$, $t_{13}=2-\delta^2$, then $t_{124}=0$, $t_{123}=3t-t^3$, and the definitive relations hold.
  \item[\rm(iii)] If $t_{13}=t_{24}=2+\delta$, then $s_{234}=s_{123}$, and
        $$s_{123}=\frac{\delta t_{13}-1}{t(\delta-1)}+\frac{t}{2}(\delta-t_{13})=\frac{\delta+1}{t(\delta-1)}\ne0.$$
        On the other hand, (\ref{eq:typeII-4}) becomes $0=2\delta s_{123}$, forcing $\delta=0$, and then it is easy to find $\det S_{123}=0$.
        This contradicts $2s_{123}^2+\det S_{123}=0$.
  \item[\rm(iv)] If $t_{13}=t_{24}=2-\delta$, then $s_{234}=s_{123}$, and
        $$s_{123}=\frac{\delta t_{13}-1}{t(\delta-1)}+\frac{t}{2}(\delta-t_{13})=\frac{\delta^2-1}{t};$$
        on the other hand, we can compute $\det S_{123}=2(\delta-1)(\delta-\delta^2)$.
        So
        \begin{align*}
        2s_{123}^2+\det S_{123}=\frac{2(\delta^2-1)^2}{\delta+2}+2(\delta-1)(\delta-\delta^2)=\frac{2(\delta-1)^2}{\delta+2}\ne 0.
        \end{align*}
\end{enumerate}

To summarize, there are seven possibilities:
\begin{enumerate}
  \item $t^2=1$, $\mathbf{x}_3=\mathbf{x}_1$, $\mathbf{x}_4=\mathbf{x}_2$;
  \item $t^2=1$, $t_{13}=t_{24}=1$, $t_{123}=t_{134}=2t$, $t_{124}=t_{234}=0$;
  \item $t^2=1$, $t_{13}=t_{24}=1$, $t_{123}=t_{134}=0$, $t_{124}=t_{234}=2t$;
  \item $t=0$, $t_{13}=-2$, $t_{24}=1$, $t_{123}=t_{124}=t_{134}=t_{234}=0$;
  \item $t=0$, $t_{13}=1$, $t_{24}=-2$, $t_{123}=t_{124}=t_{134}=t_{234}=0$;
  \item $t^2\ne0,1,3$, $t_{13}=1$, $t_{24}=4t^2-2-t^4$, $t_{123}=t_{134}=0$, $t_{124}=t_{234}=3t-t^3$;
  \item $t^2\ne0,1,3$, $t_{24}=1$, $t_{13}=4t^2-2-t^4$, $t_{123}=t_{134}=3t-t^3$, $t_{124}=t_{234}=0$.
\end{enumerate}

\begin{rmk}
\rm Case 1 can be incorporated with $\mathcal{X}'_3$ given in Section 4.1.3.

Case 3 and 5 can be incorporated with Case 6.
Case 2 and 4 can be incorporated with Case 7.
\end{rmk}

Thus, we find two more components:
\begin{align*}
\mathcal{X}_{6,1}=\{\vec{\mathsf{t}}\colon &t_{12}=t_{23}=t_{34}=t_{14}=t_{13}=1,\ t_{24}=4t^2-2-t^4, \\
&t_{123}=t_{134}=0, \ t_{124}=t_{234}=3t-t^3, \ t^2\ne 3\},   \\
\mathcal{X}_{6,2}=\{\vec{\mathsf{t}}\colon &t_{12}=t_{23}=t_{34}=t_{14}=t_{24}=1, \ t_{13}=4t^2-2-t^4, \\
&t_{124}=t_{234}=0, \ t_{123}=t_{134}=3t-t^3, \ t^2\ne 3\}.
\end{align*}

\section{Conclusions and discussions}

\subsection{The main result}

Denote a general element of $\mathcal{X}^{\rm irr}(T)$ by
$$(t_{12},t_{23},t_{34},t_{14};t_{13},t_{24};t_{123},t_{124},t_{134},t_{234}),$$
with the understanding that $t_1=t$.

Recall $\alpha=t^2-1$, and
\begin{align*}
\tau&=\big(\alpha^2-3\alpha+\frac{4}{\alpha}\big)\mu-\alpha^3+4\alpha^2-7,  \\
\sigma_+&=\big(\frac{1}{\alpha}+1\big)\mu-2+\frac{\mu}{\alpha}x, \qquad
\sigma_-=\frac{\mu^2}{\alpha}+\big(\frac{1}{\alpha}+1\big)\mu-2-\frac{\mu}{\alpha}x.
\end{align*}

\begin{thm}\label{thm:main}
The irreducible ${\rm SL}(2,\mathbb{C})$-character variety of $T$ decomposes as
$$\mathcal{X}^{\rm irr}(T)=\mathcal{X}_{1,1}\sqcup\mathcal{X}_{1,2}\sqcup\mathcal{X}_{2,1}\sqcup\mathcal{X}_{2,2}\sqcup\mathcal{X}_3
\sqcup\mathcal{X}_4\sqcup\mathcal{X}_{5,1}\sqcup\mathcal{X}_{5,2}\sqcup\mathcal{X}_{6,1}\sqcup\mathcal{X}_{6,2},$$
with
\begin{align*}
\mathcal{X}_{1,1}=\big\{&(x+1,\mu+1-x,\mu+1-x,x+1;\ 1,\tau;\ t\mu/\alpha,t\sigma_+,t\mu/\alpha,t\sigma_-)\colon \\
&\ x^2-\mu x+2\mu/\alpha+\alpha-4=0, \ \  \mu^2+(\alpha^2-4\alpha)\mu+4\alpha^2-\alpha^3=0,  \\
&\ \ \ \alpha\ne 0,2,4, \ \ \ (\alpha,\mu)\ne(2+2\sqrt{2},2\sqrt{2}),\ (2-2\sqrt{2},-2\sqrt{2})\big\},   \\
\mathcal{X}_{1,2}=\big\{&(x+1,x+1,\mu+1-x,\mu+1-x;\ \tau,1;\ t\sigma_+,t\mu/\alpha,t\sigma_-,t\mu/\alpha)\colon \\
&\ x^2-\mu x+2\mu/\alpha+\alpha-4=0, \ \  \mu^2+(\alpha^2-4\alpha)\mu+4\alpha^2-\alpha^3=0,  \\
&\ \ \ \alpha\ne 0,2,4, \ \ \  (\alpha,\mu)\ne(2+2\sqrt{2},2\sqrt{2}),\ (2-2\sqrt{2},-2\sqrt{2})\big\},   \\
\mathcal{X}_{2,1}=\big\{&(\gamma+\alpha-1,z+1,z^2/\gamma+\alpha-1,z+1;\ 1,1; \\
&\ \ \ \ t+tz\gamma/\alpha,\ t+tz\gamma/\alpha,\ t+tz^3/(\alpha\gamma),\ t+tz^3/(\alpha\gamma))\colon   \\
&\ \ \ \ z^2-\alpha z+\alpha=0, \ \  \gamma^2+((\alpha-2)z^2-2z)\gamma+z^2=0, \\
&\ \ \ \ \alpha\ne 0,2,4, \ \ \  (\alpha,z)\ne(2+2\sqrt{2},\sqrt{2}),\ (2-2\sqrt{2},-\sqrt{2})\big\},   \\
\mathcal{X}_{2,2}=\big\{&(z+1,\gamma+\alpha-1,z+1,z^2/\gamma+\alpha-1;\ 1,1; \\
&\ \ \ \ t+tz^3/(\alpha\gamma),\ t+tz^3/(\alpha\gamma),\ t+tz\gamma/\alpha,\ t+tz\gamma/\alpha)\colon  \\
&\ \ \ \ z^2-\alpha z+\alpha=0, \ \  \gamma^2+((\alpha-2)z^2-2z)\gamma+z^2=0, \\
&\ \ \ \ \alpha\ne 0,2,4, \ \ \  (\alpha,z)\ne(2+2\sqrt{2},\sqrt{2}),\ (2-2\sqrt{2},-\sqrt{2})\big\},  \\
\mathcal{X}_3=\big\{&(u,u,u,u;\ t^2-2,t^2-2;\ tu-t,tu-t,tu-t,tu-t)\colon \\
&\ \ \ \ \ u^2-(t^2+1)u+2t^2-1=0, \ t^2\ne 5\big\},   \\
\mathcal{X}_4=\big\{&(u,u,u,u;\ 2u+2-t^2,2u+2-t^2;\ \eta,\eta,\eta,\eta)\colon \\
&\ \ \ \ \ t\eta=u^2+u-t^2+1, \ \ u\ne 1, \\
&\ \ \ \ \ \eta^2+t(t^2+4-2u)\eta-5t^2u+t^4+4t^2-6=0\big\},  \\
\mathcal{X}_{5,1}=\big\{&(1,t^2-2,1,t^2-2;\ 1,1;\ 0,0,0,0)\colon t^2\ne 3\big\},  \\
\mathcal{X}_{5,2}=\big\{&(t^2-2,1,t^2-2,1;\ 1,1;\ 0,0,0,0)\colon t^2\ne 3\big\},  \\
\mathcal{X}_{6,1}=\big\{&(1,1,1,1;\ 1,4t^2-2-t^4;\ 0,3t-t^3,0,3t-t^3)\colon t^2\ne 3\big\},  \\
\mathcal{X}_{6,2}=\big\{&(1,1,1,1;\ 4t^2-2-t^4,1;\ 3t-t^3,0,3t-t^3,0)\colon t^2\ne 3\big\}.
\end{align*}
\end{thm}

\begin{rmk}
\rm Note that $t^2=3$ only occurs in $\mathcal{X}_3$ and $\mathcal{X}_4$, and $t^2=5$ only occurs in $\mathcal{X}_4$, $\mathcal{X}_{5,1}$, $\mathcal{X}_{5,2}$, $\mathcal{X}_{6,1}$, $\mathcal{X}_{6,2}$.
As an interesting consequence, if $\rho:\pi(T)\to G$ is irreducible with $\chi_\rho\notin\mathcal{X}_4$ and $t^2=3$ (resp. $t^2=5$), then $\rho$ must factor through $\pi(4_1)$ (resp. $\pi(3_1)$); see the next remark.

The intersection of any two of the closures of $\mathcal{X}_4$, $\mathcal{X}_{5,1}$, $\mathcal{X}_{5,2}$, $\mathcal{X}_{6,1}$, $\mathcal{X}_{6,2}$
consists of $(1,1,1,1;1,1;0,0,0,0)$ for $t=\pm\sqrt{3}$, which are reducible.

The map $f^\ast$ interchanges $\mathcal{X}_{i,1}$ with $\mathcal{X}_{i,2}$ for $i=1,2,5,6$, and fixes $\mathcal{X}_3$, $\mathcal{X}_4$.
\end{rmk}

\begin{rmk}[Compared with \cite{PP20} Section 6]
\rm As we have checked, the closure of $\mathcal{X}_3$ is the ``figure eight knot component". Indeed, $\mathcal{X}_3$ consists of $\chi_\rho$ for $\rho$'s factoring through the homomorphism
\begin{align*}
\pi(T)\twoheadrightarrow\pi(4_1)&=\langle s_1,s_2\mid s_1s_2s_1^{-1}s_2s_1=s_2s_1s_2^{-1}s_1s_2\rangle,  \\
&x_1,x_3\mapsto s_1, \qquad  x_2,x_4\mapsto s_2.
\end{align*}

We also have checked that the closures of $\mathcal{X}_{5,1}$, $\mathcal{X}_{5,2}$, $\mathcal{X}_{6,1}$, $\mathcal{X}_{6,2}$ are the ``trefoil-like" components respectively given by (5), (7), (6), (4) in \cite{PP20}.

Furthermore, the closures of $\mathcal{X}_{1,1}$, $\mathcal{X}_{1,2}$, $\mathcal{X}_{2,1}$, $\mathcal{X}_{2,2}$ are called the ``four more components", equivalent under the action of the symmetry group of $T$.
\end{rmk}

Due to their complexities and mysteries, the four components $\mathcal{X}_{1,1}$, $\mathcal{X}_{1,2}$, $\mathcal{X}_{2,1}$, $\mathcal{X}_{2,2}$ may be called ``exotic components".

We hope to answer the following question in the future:
\begin{que}
Do $\mathcal{X}_{1,1}$, $\mathcal{X}_{1,2}$, $\mathcal{X}_{2,1}$, $\mathcal{X}_{2,2}$ encode any information on geometric structures of the exterior of $T$?
\end{que}

\subsection{The excellent component}

The closure of $\mathcal{X}_4$ is the excellent component $\mathcal{X}_\star$ given in \cite{HLM00}.
To be precise, when $t\ne 0$, we can eliminate $\eta$ to obtain
$$u^4+(2-2t^2)u^3+(t^4+3)u^2+(2-2t^4)u+2t^4-4t^2+1=0,$$
whose left-hand-side coincides with $-\hat{r}_{E[8_{18}]}(t,u)$ given on \cite[Page 9]{HLM00}.

Note that the ``canonical component" obtained in \cite[Page 20]{PP20}, as claimed by the authors, is a double branched cover of $\mathcal{X}_\star$.

When $t^2=3+2\epsilon\sqrt{2}$ with $\epsilon\in\{\pm1\}$,
the point given by $u=1+\epsilon\sqrt{2}$ and $\eta=(u+1)/t$ is very special, in the sense that it is a limit point of each of
$\mathcal{X}_{1,1}$, $\mathcal{X}_{1,2}$, $\mathcal{X}_{2,1}$, $\mathcal{X}_{2,2}$.
Furthermore, the point given by $t=u=1-\sqrt{2}$, $\eta=-\sqrt{2}$ is the limit of hyperbolicity described on \cite[Page 9]{HLM00}.

Also contained in $\mathcal{X}_4$ are some characters of number-theoretic interests:
\begin{itemize}
  \item when $t=0$, $u=(-1+\epsilon\sqrt{3}i)/2$, $\eta=\epsilon'\sqrt{6}$, with $\epsilon,\epsilon'\in\{\pm1\}$;
  \item when $t^2=1$, $u^4+4u^2-1=0$, i.e. $u^2=-2\pm\sqrt{5}$;
  \item when $t^2=3$, since $u=1$ is excluded, one has $u=1\pm\sqrt{6}i$;
  \item when $t^2=5$, $u^4-8u^3+28u^2-48u+31=0$, whose roots are
        $$2\pm\sqrt{\sqrt{5}-2}, \qquad  2\pm i\sqrt{\sqrt{5}+2};$$
  \item when $t^4-6t^2+4=0$, i.e. $t^2=3\pm\sqrt{5}$, one has $t_{12}=t^2/2$, $\mathbf{x}_3=\mathbf{x}_1^{-1}$,
        and $\mathbf{x}_4=\mathbf{x}_2^{-1}$.
\end{itemize}
It will be interesting to find geometric interpretation for any of them.

\subsection{Parabolic representations}

We are going to find equivalence classes of parabolic representations, meaning those with $t=\pm2$ up to conjugacy and multiplication by $-\mathbf{e}$. We may just set $t=2$.
Let $\mathcal{X}^{\rm par}_{1,1}=\{\chi\in\mathcal{X}_{1,1}\colon t=2\}$, and so forth.
\begin{itemize}
  \item Each of $\mathcal{X}^{\rm par}_{1,1}$, $\mathcal{X}^{\rm par}_{1,2}$ consists of $4$ elements, determined by $\mu,x$, with
        $$\mu=\frac{1}{2}(3+3\epsilon\sqrt{3}i), \quad  x=\frac{1}{4}\big(3+3\epsilon\sqrt{3}i\pm\sqrt{2\epsilon\sqrt{3}i-18}\big),
        \quad \epsilon\in\{\pm1\}.$$
  \item Each of $\mathcal{X}^{\rm par}_{2,1}$, $\mathcal{X}^{\rm par}_{2,2}$ consists of $4$ elements, determined by $z,\gamma$, with
        $$z=(3\pm\sqrt{3}i)/2,  \quad   \gamma=\frac{1}{4}\big(3+\epsilon\sqrt{3}i\pm\sqrt{30\epsilon\sqrt{3}i-18}\big),
        \quad \epsilon\in\{\pm1\}.$$
  \item $\mathcal{X}^{\rm par}_3$ consists of $2$ elements, determined by $u=(5\pm \sqrt{3}i)/2$.
  \item $\mathcal{X}^{\rm par}_4$ consists of $4$ elements, determined by the roots of
        $$u^4-6u^3+19u^2-30u+17=0,$$
        namely,
        $$\frac{1}{2}\Big(3\pm\sqrt{8\sqrt{2}-11}\Big), \qquad     \frac{1}{2}\Big(3\pm i\sqrt{8\sqrt{2}+11}\Big).$$
  \item Each of $\mathcal{X}^{\rm par}_{5,1}$, $\mathcal{X}^{\rm par}_{5,2}$, $\mathcal{X}^{\rm par}_{6,1}$,
        $\mathcal{X}^{\rm par}_{6,2}$ consists of a single element.
\end{itemize}

Thus, there are 26 equivalence classes of parabolic representations in total. This is consistent with \cite[Proposition 8.1]{CKKY22}.

\bigskip

\noindent
Haimiao Chen (orcid: 0000-0001-8194-1264)\ \ \ \ {\it chenhm@math.pku.edu.cn} \\
Department of Mathematics, Beijing Technology and Business University, \\
Liangxiang Higher Education Park, Fangshan District, Beijing, China.

\end{document}